\documentclass[A4paper,11pt]{amsart}
\usepackage[english] {babel}
\usepackage{amsmath,amssymb,amscd}
\usepackage{amsthm}
\usepackage{subfig}
\usepackage{graphicx}
\usepackage{xcolor}
\usepackage{hyperref}
\usepackage[utf8]{inputenc}
\usepackage{fancybox}
\usepackage{caption}
\usepackage[all]{xy}

\usepackage{bookmark,hyperref}
\hypersetup{
     colorlinks=true,
     linkcolor=blue,
     filecolor=blue,
     citecolor = blue,
     urlcolor=cyan,
     }

\newcommand{\A}{{\mathbb A}}
\newcommand{\N}{{\mathbb N}}
\newcommand{\cH}{{\mathbb H}}
\newcommand{\calH}{{\mathcal H}}
\newcommand{\Z}{{\mathbb Z}}
\newcommand{\R}{{\mathbb R}}
\newcommand{\C}{{\mathbb C}}
\newcommand{\bs}{{\mathbb S}}
\newcommand{\Q}{{\mathbb Q}}
\newcommand{\bO}{{\mathbf\Omega}}
\newcommand{\CQ}{{\C^*_\Q}}
\newcommand{\SQ}{{S^{1}_{\Q}}}
\newcommand{\PQ}{\C P^{1}_\Q}
\newcommand{\bS}{\mathbb{S}}
\newcommand{\bu}{{\mathbf{u}}}
\newcommand{\bv}{{\mathbf{v}}}
\newcommand{\bw}{{\mathbf{w}}}
\newcommand{\bd}{{\mathbf{d}}}
\newcommand{\bz}{{\mathbf{z}}}
\newcommand{\I}{{\mathcal O}}
\newcommand{\La}{{\boldsymbol\Lambda}}
\newcommand{\gln}{{\text{GL}(n,\C)}}
\newcommand{\Lo}{{{\boldsymbol\Lambda}_0^\infty(G)}}
\newcommand{\Ld}{{{\boldsymbol\Lambda}^\infty(G)}}

\def\hmath$#1${\texorpdfstring{{\rmfamily\textit{#1}}}{#1}}

\newcommand{\dem}{{\em Proof: \;}}
\newcommand{\fdem}{\hfill $\square$}

\newtheorem{definition}{Definition}
\newtheorem{example}{Example}

\newtheorem{remark}{Remark}

\newtheorem{theorem}{Theorem}
\newtheorem{corollary}{Corollary}
\newtheorem{proposition}{Proposition}
\newtheorem{lemma}{Lemma}

\def\noi{\noindent}
\def\ie{\emph{i.e.,\,}}

\title[Adelic toric varieties and adelic loop groups]{Adelic toric varieties and adelic loop groups}

\author{J. M. Burgos  and  A. Verjovsky}



\begin{document}

\maketitle
{\footnotesize
\centerline{\it To Dennis Sullivan on occasion of his 80th  birthday and for teaching us the dictum:} 
\centerline{\it  ``solenoidal manifolds are diffuse versions of manifolds''}

\begin{abstract}
We study the proalgebraic space which is the inverse limit of all finite branched covers over a normal toric variety with branching set the invariant divisor under the algebraic torus action. These are completions (compactifications) of the adelic abelian algebraic group which is the profinite completion of the algebraic torus. 
We prove that the vector bundle category of the proalgebraic toric completion of a toric variety is the direct limit of the respective categories of the finite toric varieties coverings defining the completion. In the case of the complex projective line
we obtain as proalgebraic completion the adelic projective line $\PQ$. We define 
holomorphic vector bundles over $\PQ$. We also introduce
the smooth, Sobolev and Wiener adelic loop groups and the corresponding Grassmannans; we describe their properties and prove Birkhoff's factorization for these groups.
We prove that the adelic Picard group of holomorphic line bundles is isomorphic to the rationals and prove the Birkhoff-Grothendieck splitting theorem for holomorphic bundles of higher rank over $\PQ$.

\end{abstract}

\tableofcontents

\section{Introduction}\label{proalgebraic}

A \emph{proalgebraic variety} is a projective system of complex algebraic varieties
$$f_{ij}:V_j \to  {V_i} \quad i, j\in{I}, \,\,i\leq{j} $$
where $(I,\leq)$ is a directed poset. Each $V_i$ is a complex algebraic variety and the bonding maps $f_{ij}, \, i\leq{j}$ are algebraic maps. 

Proalgebraic varieties play an important role in several branches of mathematics in particular in the theory of dynamical systems. For instance, in \cite{Gr1} M. Gromov studies the \emph{surjunctivity} of algebraic maps $f:V \to V$ (\ie regular algebraic maps which are injective are necessarily also surjective) and in \cite{Gr2}, using the inverse system
$$V\overset{f}\longleftarrow{V}\overset{f}\longleftarrow{V}\cdots$$
Gromov initiated what he calls \emph{symbolic algebraic geometry} because one can study the properties of $f$ by means of the associated shift map in the inverse system \ie  through symbolic dynamics.

Examples of these profinite structures appear constantly in dynamical systems, for instance, the dyadic solenoid of Van-Danzig Vietoris is an invariant set of a diffeomorphism $f:{}\mathbb S^3\to{\mathbb S}^3$ which is the Smale--Williams expansive attractor of an Axiom A diffeomorphism of the 3-sphere.

\emph{Riemann surface laminations} which are transversally like the Cantor set were introduced into conformal dynamics by Sullivan \cite{Su1}. These Riemann surface laminations play a role in the dynamics similar to the role played by Riemann surfaces in the action of Kleinian groups (\cite{Su2},  \cite{Su3}). 
The type of laminations we consider in this paper are called  \emph{solenoidal manifolds}  (\cite{Su4,Ve}). These are Hausdorff topological spaces $X$ equipped with a covering $\{U_i\}_{i\in\mathcal{I}}$ by open sets
and coordinate charts $\phi_i:U_i\to{{D_i}\times{T_i}}$, where $D_i$ is homeomorphic to a domain in $\R^n$ and $T_i$ is homeomorphic 
to the Cantor set (in our case a Cantor group).
The transition maps $\phi_i\circ{\phi_j}^{-1}=\phi_{ij}:\phi_j(U_i\cap{U_j})\to\phi_i(U_i\cap{U_j})$
are required to be homeomorphisms that take plaques to plaques, \ie they are of the form:
$$
\phi_{ij}(x,t)=(f_{ij}(x,t),\gamma_{ij}(t)).
$$
For a fixed $t=t_0$ one may impose some regularity on the maps $\phi_{ij}(\cdot\,,t_0)=(f_{ij}(\cdot\,,t_0),\gamma_{ij}(t_0))$ in their domain of definition which can be identified with a map between two open sets in $\R^n$. These are objects which locally look like a product of open $n$-disk with a Cantor set. When $n=2$, identifying $\R^2$ with $\C$ and requiring that the maps between plaques to be holomorphic, we obtain a Riemann surface lamination. 

Sullivan associated a Riemann surface lamination to any $C^2$-smooth expanding map of the circle. M. Lyubich and Yair Minsky \cite{LM} used laminations in dimension 3 called \emph{hyperbolic orbifold 3-laminations} (a 3-dimensional solenoidal manifold) to study holomorphic dynamics. Both Dennis Sullivan and Misha Lyubich-Yair Minsky obtained remarkable results in conformal dynamics using these solenoidal manifolds.
 
 The first four sections of this paper develop and study the notions of 
 {\bf adelic toric variety}. These varieties are completions (compactifications) of proalgebraic, abelian groups, which are defined and commented in section \ref{abelian-pro-Lie}. These groups are abelian pro-Lie groups. 
 
 \medskip
  In section \ref{adelic_solenoid} we define the {\bf adelic solenoid} $\SQ$. This is the compact abelian group which is the Pontryagin dual of the additive group of the rationals 
 $\Q,+$, with the discrete topology. Is a compact solenoid of dimension one. It is a ``diffuse''
 version of the circle $\bS^1$. 
 {\it This metaphor is at the basis of our definition of the 
 adelic loop group of a Lie group}.
 
  \medskip
 In section \ref{proalgebraic_completion} we define proalgebraic normal toric varieties 
and proalgebraic completion of a toric variety.
 
 If $X^\mathcal{F}$ is a \emph{normal} complex toric projective variety of complex dimension $n$ determined by the fan 
 $\mathcal{F}$ then, by definition, there is 
 a holomorphic action of the algebraic torus $(\C^*)^n=\C^*\times\cdots\times\C^*\,$ 
 (see Definition (\ref{profinite-algebraic-torus}) below) with a Zariski-dense principal orbit $U\subset{X^\mathcal{F}}$. The
 complement $D:=X^\mathcal{F}-U$ is a divisor which is invariant under the action. The profinite 
 completion $X^{\mathcal{F}}_\Q$ of $X^\mathcal{F}$ is the inverse limit of the tower of finite coverings branched over $D$. Normality
 implies that the action of $(\C^*)^n$ lifts to any such finite  branched covering and therefore these coverings are
 also normal toric complex varieties \cite{AP}. 
 
 The actions on the covers determine an action of 
 the proalgebraic (or adelic) torus (see Section (\ref{abelian-pro-Lie}))
 ${(\C^*_{\Q})}^n=\C^*_{\Q}\times\cdots\times\C^*_{\Q}$ on $X^{\mathcal{F}}_\Q$, where $\C^*_\Q$ is the profinite completion
 of $\C^*$ (see Section (\ref{adelic_solenoid})). Thus the proalgebraic toric varieties are equivariant
 compactifications of the ``solenoidal'' abelian group ${(\C^*_\Q)}^n$.

Consider the index set $I$ of the projective system of finite coverings $X^{\mathcal{F}}_i \rightarrow X^{\mathcal{F}}$ branched over $D$. By definition, for every index $i$ in $I$, there is a canonical map $\pi_i: X^{\mathcal{F}}_\Q\rightarrow X^{\mathcal{F}}_i$. In subsection (\ref{vector_section}) we prove the following

\begin{theorem}\label{main_vector}
For every index $i$ in $I$ and every natural $r$, the pullback by $\pi_i$ is a full and faithful functor between the rank $r$ vector bundle categories
$$\pi_i^{*}: \mbox{Vect}_r\left(X^{\mathcal{F}}_i\right)\rightarrow \mbox{Vect}_r\left(X^{\mathcal{F}}_\Q\right)$$
and the collection of these functors define an inductive system whose limit is equivalent to the rank $r$ vector bundle category of the proalgebraic completion $X^{\mathcal{F}}_\Q$, that is to say
$$\mbox{Vect}_r\left(X^{\mathcal{F}}_\Q\right)\cong \varinjlim_{i\in I} \mbox{Vect}_r\left(X^{\mathcal{F}}_i\right).$$
\end{theorem}

The following corollaries are immediate consequences of Theorem \ref{main_vector}.

\begin{corollary}
For every index $i$ in $I$, the pullback by $\pi_i$ induce a morphism of $\mbox{K}$--rings $\mbox{K}(X^{\mathcal{F}}_i)\rightarrow \mbox{K}(X^{\mathcal{F}}_\Q)$ and the collection of them define an inductive system whose limit is isomorphic to $\mbox{K}(X^{\mathcal{F}}_\Q)$.
\end{corollary}

The homology, cohomology and characteristic classes of toric varieties have been studied extensively, see for instance \cite{Yokura1, Yokura2}, \cite{BBFK1, BBFK2}. 

The following result reduces the problem of calculating the Chern class of a complex vector bundle over the proalgebraic completion to doing so on a finite branched covering.

\begin{corollary}
Consider a rank $r$ complex vector bundle $E\rightarrow X^{\mathcal{F}}_\Q$. Then, there is an index $l$ and a rank $r$ complex vector bundle $E_l\rightarrow X^{\mathcal{F}}_l$ such that $c(E)= \pi^*_l\,c(E_l)$ where $c$ denotes the total Chern class.
\end{corollary}

The inverse image of the divisor $D$ by the canonical projection $X^{\mathcal{F}}_\Q \rightarrow X^{\mathcal{F}}$ constitutes the set of singularities of the proalgebraic completion. In the case of the proalgebraic completion of a projective toric variety, these singularities are combinatorially organized by the incidence relations of the faces of the Delzant polytope of $X^{\mathcal{F}}$, see (section (\ref{proalgebraic_completion}), Proposition \ref{Delzant_prop}). 

\medskip
We wonder whether Theorem \ref{main_vector} admits a proof depending only on the classifying space or universal object and holds on other categories other than the vector bundle one. In view of this question, we believe that Theorem \ref{main_vector} does not constitute its full generality.

\medskip
Section \ref{hol-vect-bun} deals with the notions and properties of holomorphic functions and vector bundles on $\PQ$. Subsection \ref{lp} introduces the notion of holomorphic functions on adelic punctured disks and annuli via an adapted version of Laurent-Puiseux series. Subsection
\ref{wn} defines the notion of {\bf winding numbers}, which are topological invariants
with values in $\Q$ and play the role of Chern numbers for adelic line bundles.

In section \ref{adelic-loop-groups} we define the adelic loop groups of a Lie group 
$G$ with emphasis on the Lie groups $\gln$ and $U(n)$. In section \ref{Wienerloops} 
we study loops with coefficients in the Wiener algebra. 
In section \ref{ksloops} we provide 
with a  Kähler structure the loop group with coefficients in Sobolev space $H^1$ and define
the Energy Functional of these groups with this Kähler metric in subsection \ref{energy}.
Using the ``descending'' semiflow of the gradient of the energy functional we prove a Birkhoff decomposition of the adelic group loop in subsection \ref{bd}. 
In subsection \ref{filtration} we describe a filtration of the continuous loop group
according to different degrees of regularity.
These groups are interesting in their own right but are used in this paper to prove Birkhoff factorization theorem and Iwasawa factorization in section \ref{factorizations}.
In subsection \ref{sobgrass} we define the Sobolev Grassmannian using the fact, developed in 
section \ref{loop-operators}, that loops can be regarded as operators of certain types
in a Sobolev space. This implies the remarkable fact, valid for the classic based loop space of the unitary group,
that the Grassmannian is biholomorphic to the Kähler manifold of based Sobolev loops.
on the unitary group $U(n)$. In section \ref{BF} we prove the Brirkhoff factorization theorem for adelic loops and use it to prove in section \ref{BG} the Birkhoff-Grothendieck decomposition theorem: any holomorhic vector bundle over $\PQ$ of rank $n$ splits as a direct sum of $n$ holomorphic line bundles. The splitting is unique up to reordering the summands.

\noi The last section {\it Concluding Remarks} (\ref{concrmks}) indicates some possible paths for future research on the subject.
\label{bd}

\medskip
Finally, we would like to point that the {\bf ``Adelic''} adjective in the title
refers to the fact that the proalgebraic toric completions in this paper are completions
of the abelian adelic algebraic group corresponding to 
$(GL(1,\C))^n=(\C^*)^n=\C^*\times\cdots\times\C^*$ (\cite{Bo, RV, T, We}). Its maximal compact subgroup is $(\SQ)^n$, and $\SQ$ is the {\bf adèle class group} of the rationals. 

\section{Adelic solenoid}\label{adelic_solenoid}

For every $n,m\in \Z^+$ such that $n$ divides $m$, then there exists a covering map $p_{n,m}:S^1\to S^1$ such that $p_{n,m}(z)=z^{m/n}$. 

This determines a projective system of covering spaces $\{S^1,p_{n,m}\}_{n,m \geq 1, n|m}$ whose projective limit is the \textsf{universal one--dimensional solenoid} or \textsf{adelic solenoid}:
\[ S^1_\Q:=\lim_{\underset{p_{n,m}}\longleftarrow}S^1. \] 
Thus, as a set, $S^1_\Q$ consists of sequences $\left( z_n\right)_{n\in\N,\, z\in{S^1}}$ such that $p_{n,m}(z_m)=z_n$ if $n$ divides $m$.

The canonical projections of the inverse limit are the functions 
$S^1_\Q\overset{\pi_n}\to S^1$ defined by 
\[
\pi_n\left(\left(z_j\right)_{j\in\N}\right)=z_n 
\]

\noi and they define the solenoid topology as the initial topology of the family. The solenoid is an abelian topological group and each $\pi_n$ is an epimorphism.  In particular each $\pi_n$ is a character which determines a locally trivial $\hat{\Z}$--bundle structure where the group

\[
\hat{\Z}:=\lim_{\underset{n|m}\longleftarrow}
\left\{\Z/m\Z\,\overset{p_{n,m}}\longrightarrow\,\Z/n\Z\right\}, 
\]

\noi (where $p_{n,m}:\Z/m\Z\longrightarrow\,\Z/n\Z$ is the canonical epimorphism when $n|m$), is the profinite completion of $\Z$. Then, $\hat{\Z}$ is a compact, perfect and totally disconnected abelian topological group homeomorphic to the Cantor set. 

Since $\hat{\Z}$ the  profinite completion of $\Z$, it admits a canonical inclusion of $\Z\subset\hat\Z$ whose image is dense.
We have an inclusion $\hat\Z\overset{\phi}\to{S^1_\Q}$ and a short exact sequence:
$$0\to{\hat\Z}\overset{\phi}\rightarrow S^1_\Q\overset{\pi_1} \rightarrow S^1\to1$$

The solenoid $S^1_\Q$ can also be realized as the orbit space of the $\Q$--bundle structure $\Q \hookrightarrow \mathbb{A} \to \A/\Q$, where $\A$ is the {\bf adele group} of the rational numbers which is a locally compact Abelian group, $\Q$ is a discrete subgroup of $\A$ and $\A/\Q \cong S^1_\Q$ is a compact Abelian group (\cite{RV}). From this perspective, $\A/\Q$ can be seen as a projective limit whose $n$--th component corresponds to the unique covering of degree $n\geq 1$ of $S^1_\Q$.

By considering the properly discontinuously diagonal free action of $\Z$ on $\hat{\Z}\times\R$  given by
\[ n\cdot(x,t)=(x+n,t-n), \quad (n\in \Z, \, x\in\hat\Z, \,t\in\R)\]
the solenoid $S^1_\Q$ is identified with the orbit space $\hat{\Z}\times_{\Z} \R$. Here, $\Z$ is acting on $\R$ by covering transformations and on $\hat{\Z}$ by translations. The path--connected component of the identity element $1\in S^1_\Q$ will be called the \textsf{baseleaf} \cite{Od} and it is a densely immersed copy of $\R$.

Hence $S^1_\Q$ is a compact, connected, abelian topological group and also a one-dimensional lamination where each ``leaf" is a simply connected one-dimensional manifold, homeomorphic to the universal covering space $\R$ of $S^1$, and a typical ``transversal" is isomorphic to the Cantor group $\hat{\Z}$. The solenoid $S^1_\Q$ also has a leafwise $\mathrm{C}^\infty$ Riemannian metric (i.e., $\mathrm{C}^\infty$ along the leaves) which renders each leaf isometric to the real line with its standard metric $dx$. So, it makes sense to speak of a rigid translation along the leaves. The leaves also have a natural order equivalent to the order of the real line hence also an orientation.

Summarizing the above discussion we have the commutative diagram: 
\begin{equation}\label{diagram_I}
\xymatrix{
    S^{1}_{\Q}= \varprojlim S^{1} \quad \ldots\ar[r] &	 S^{1} \ar[r]^{p_{m,n}}		& 	S^{1} \ldots\ar[r] 		& 	 S^{1}   \\
	\hat{\Z}= \varprojlim\Z/n\Z \quad \ar@{^{(}->}[u]^{\phi} \ldots\ar[r]  & 	 \Z/n\Z \ar@{^{(}->}[u]_ {l\,\mapsto{e^{2\pi il/n}}} \ar[r]^{p_{m,n}} 	&	\Z/m\Z \ar@{^{(}->}[u]_{l\,\mapsto{e^{2\pi il/m}}} \ldots\ar[r]	&	  \{0\} \ar@{^{(}->}[u]_{0\,\mapsto1}}
\end{equation}
where $\hat{\Z}$ is the adelic profinite completion of the integers and the image of the group monomorphism $\phi:(\hat{\Z},+)\rightarrow (S^{1}_{\Q},\cdot)$ is the \textsf{principal fiber}. We notice that $\pi_{n}(x)= \pi_{n}(y)$ implies $\pi_{n}(y^{-1}x)=1$ and therefore $y^{-1}x= \phi(a)$ where $a\in n\hat{\Z}$ for some $n\in\Z\subset\hat{\Z}$.

We define the \textsf{baseleaf} as the image of the monomorphism $\nu:\R \rightarrow S^{1}_{\Q}$ defined as follows:
\begin{equation}\label{diagram_II}
\xymatrix{
    S^{1}_{\Q}= \varprojlim S^{1} \quad \ldots\ar[r] &	 S^{1} \ar[r]^{p_{m,n}} 	
    & 	S^{1} \ldots\ar[r] 		& 	 S^{1}   \\
	\R \ar@{^{(}->}[u]_{\nu} \quad \ldots\ar[r]^{=}  & 	 \R \ar[u]_{t\,\mapsto{e^{it/n}}} \ar[r]^{=} 	&	\R \ar[u]_{t\,\mapsto{e^{it/m}}} \ldots\ar[r]^{=}	&	 \R \ar[u]_{t\,\mapsto{e^{it}}}}.
\end{equation}
In particular, the immersion $\nu$ is a group morphism and comparing the diagrams \eqref{diagram_I} and \eqref{diagram_II}, we have $\nu(2\pi n)= \phi(n)$ for every integer $n$.

\begin{definition}[Adelic exponential and canonical flow]\label{exp}
Let $\boldsymbol{Exp}:\R\times\hat{\Z}\rightarrow S^{1}_{\Q}$ such that 
 $\boldsymbol{Exp}(t,a)= \nu(t)\phi(a)$. Then $\boldsymbol{Exp}$ is an epimorphism with kernel
the Cantor group $\hat\Z$. This is the adelic version of the map 
$t\mapsto{e^{2\pi\boldsymbol{i}{t}}}$.
Consider the translation flow $\hat\varphi_t$ on $\R\times\hat{\Z}$
given by $\hat\varphi_t(s,\bz)=(s+t,\bz)$. This flow commutes with the
map $(t,\bz)\mapsto(t+1,\bz+\mathbf1)$ so it defines a translation flow 
$\varphi_t:\SQ\to\SQ$ which preserves the 1-dimensional leaves of the lamination
structure of $\SQ$; it amounts to translations by elements in the base leaf. This flow is called {\bf the canonical flow} of $\SQ$.
The vector field (in the solenoidal sense) is the unit vector field along the
leaves. The flow is in fact the one-parameter group of translations by the
connected component of the identity (the base leaf). 
\end{definition}

\noi {\bf We will always treat $\SQ$ as a multiplicative abelian group.}

\begin{remark}\label{integeradeles} If $\mathbb{A} _{\mathbb{Z}}$ denotes the ring of integral adeles \cite{RV} then: 
 $\mathbb{R}\times{\hat{\mathbb {Z}}}=\mathbb{A}_{\mathbb{Z}}=
 \mathbb{R}\times\prod_{p}\mathbb{Z}_{p}$. 
 The map $\Z\overset{i}\hookrightarrow\mathbb{R}\times{\hat{\mathbb {Z}}}$, 
 $\,n\mapsto(n,\mathbf{n})$,
 where $\mathbf{n}$ corresponds to the natural inclusion of $\Z$ into $\hat\Z$,
injects $\Z$ into a discrete co-compact subgroup $\Gamma$ so that 
$\SQ=(\R\times\hat\Z)/\Gamma$. The subgroup 
$\left\{(0,\mathbf{n}): n\in\Z\right\}$ is dense in 
$\left\{0\right\}\times\hat\Z$. This implies that the canonical flow is minimal.
\end {remark}
\medskip
The previous constructions can be extended to the inverse limit of the multiplicative group $\C^{*}$ as well as to the multiplicative semigroup $\C$. In other words, if $z_n:\C\to\C$ denotes the map $z\mapsto{z^n}$ we can take the inverse limits under the partial order of divisibility of the integers of the finite coverings of $\C^*$ and the branched coverings of $\C$:   
\begin{definition}\label{C*QandCQ}
$\varprojlim\, \C^*\overset{def}=\C^{*}_\Q$ and
$\varprojlim \,\C\overset{def}=\C_\Q$, respectively.
\end{definition}

In the first case we obtain the topological abelian group: the \textsf{algebraic solenoidal group} $\C^{*}_\Q$. This group is laminated by densely immersed copies of $\C$. In the second case, we obtain the semigroup $\C_\Q$ which is a ramified covering of $\C$ at zero and topologically it is homeomorphic to the open cone over the adelic solenoid. The singularity of a cone over a solenoidal torus $S^{1}_{\Q}\times\ldots  S^{1}_{\Q}$ will be called a \textsf{solenoidal cusp} or simply a cusp. 

The following is the complex version of the diagram \ref{diagram_I}:
\begin{equation}\label{diagram_III}
\xymatrix{
    \C^{*}_{\Q}= \varprojlim \C^{*} \quad \ldots\ar[r] & \C^{*} \ar[r]^{p_{m,n}}		& 	\C^{*} \ldots\ar[r] 		& 	 \C^{*}   \\
	\hat{\Z}= \varprojlim \Z/n\Z  \quad\ar@{^{(}->}[u]^{\phi} \ldots\ar[r]  & \Z/n\Z \ar@{^{(}->}[u]_ {l\,\mapsto{e^{2\pi il/n}}} \ar[r]^{p_{m,n}} 	&	\Z/m\Z \ar@{^{(}->}[u]_{l\,\mapsto{e^{2\pi il/m}}} \ldots\ar[r]	& \{0\} \ar@{^{(}->}[u]_{0\,\mapsto1}}
\end{equation}

\section{Algebraic abelian pro-Lie groups and proalgebraic normal toric varieties.}\label{abelian-pro-Lie}

A topological group is called an \emph{abelian pro-Lie group} (\cite{HM, S}) if it is isomorphic to a closed
subgroup of a product (possibly infinite) of finite-dimensional real abelian Lie groups. This class of groups is closed
under the formation of arbitrary products and closed subgroups and forms a complete
category. In this paper we use the following specific abelian pro-Lie group:
\begin{definition}\label{profinite-algebraic-torus}
 The $n$-dimensional \emph{algebraic profinite torus} is the profinite completion
of the algebraic torus $(\C^*)^n=\C^*\times\cdots\times\C^*\,$ ($n\,$ times)
\end{definition}
\begin{definition} For any $n\in\N$ the group ${(\C^*_\Q)}^n=\C^*_\Q\times\cdots\times\C^*_\Q$ ($n$ factors) is called the proalgebraic (or adelic) torus.
\end{definition}
Let us consider the profinite completion of $(\C^*)^n$. The fundamental group of  $(\C^*)^n$ is isomorphic 
to the free abelian group of rank $n$, $\,\Z^n=\Z\times\cdots\times\Z$. Since the finite coverings of  $(\C^*)^n$ 
correspond to the subgroups (sublattices) $\Lambda$ of $\Z^n$, of finite index,  and the covering
can be chosen to be  epimorphism $p_\Lambda:(\C^*)^n\to(\C^*)^n$ of $(\C^*)^n$ onto itself the profinite completion
is an abelian group. If we chose the canonical basis $\left\{\mathbf{e}_1,\cdots,\mathbf{e}_n\right\}$
(as a free module over $\Z$) of $\Z^n$ a lattice $\Lambda$ corresponds to the lattice generated by the columns of an $(n\times{n})$-matrix $A$
with integer entries and non-vanishing determinant which is an integer $k$. By Cramer's rule the matrix $B:=kA^{-1}$ is a
non singular matrix with integer entries. Then $BA=kI_N$ and since the columns of $BA$ are integer linear combinations
of the columns of $A$ it follows that $\Lambda$ contains as sublattice the ``cubic'' lattice in $\Z^n$ generated by
 $\left\{k\mathbf{e}_1,\cdots,k\mathbf{e}_n\right\}$. Of course the lattice generated by $\left\{k\mathbf{e}_1,\cdots,k\mathbf{e}_n\right\}$ contains as a sublattice the lattice generated by $\left\{lk\mathbf{e}_1,\cdots,lk\mathbf{e}_n\right\}$ for any integer
 $l\neq0$. Thus the set of cubic lattices is cofinal in the set of all latices.
 Therefore we have
 \begin{proposition}\label{profinite-is-a-product} The profinite completion of $(\C^*)^n$ is isomorphic to the product $(\C^*_\Q)^n=\C^*_\Q\times\cdots\times\C^*_\Q$ ($n$ factors).
 \end{proposition}
Let us recall that a \emph{normal toric variety} $X$ of complex dimension $n$ is a normal projective variety which has a holomorphic action $(\C^*)\times{X}$ with an open principal orbit $U\cong(\C^*)$. It follows that $D:=X-U$ 
is a divisor which is invariant under the action. $D$ is stratified in a finite number of strata which are projective subvarieties invariant under the action, $D$ is a union of a finite set of invariant divisors $D_1,\cdots,D_k$ \cite{BBFK}.

We recall the following theorem in \cite{AP} (stated here for \emph{complex} normal toric varieties):

\begin{theorem}\label{toric-cover} Let $X$ be a normal complete toric variety defined over $\C$, of dimension $n$, and let
$f:Y\to{X}$ be a connected branched cover such that the divisorial part of the branch locus of $f$
is contained in the union of the invariant divisors $D_1,\cdots, D_k$. Then $f$ is a toric cover.
\end{theorem}
The theorem establishes that if $X$ is a normal complex toric variety with principal orbit $U\cong(\C^*)^n$, then
if $f:Y\to{X}$ is a finite branched covering with branch locus $D$ then:

\begin{enumerate}
\item $Y$ can be endowed with the structure of a projective normal toric variety with an action of $(\C^*)^n$ with principal orbit $V=f^{-1}(U)$.
\item $f$ is an algebraic equivariant map respect to the actions of $(\C^*)^n$ on $Y$ and $X$, respectively. 
\end{enumerate}

\begin{definition}  The equivariant branched covering map $f$ as in Theorem \ref{toric-cover}  
is called a \emph{toric-cover}.
\end{definition}
Toric covers of $X$ are in one-to-one correspondence with finite-index subgroups $\Lambda$ of $\pi_1(U)=\Z^n$. They form
a directed set under inclusion. Then, for each $\Lambda\subset\Z^n$ subgroup of finite index we have a toric cover 
$f:X_\Lambda\to{X}$.
Therefore we can construct the corresponding projective limit $X_\Q\underset{def}=\underset{\Lambda}\varprojlim\ X_\Lambda$.
\begin{definition}
$X_\Q$ is called the \emph{proalgebraic toric completion} of the normal toric variety $X$.
\end{definition}
\begin{remark} The proalgebraic toric completion of $X$ is not the profinte completion under \emph{every} branched cover
with branch locus a divisor. We only consider toric covers.
\end{remark}
The fact that the coverings used to construct $X_\Q$ are equivariant implies that there exists an action of the profinite 
completion of $(\C^*)^n$:

\begin{proposition} The proalgebraic toric completion $X_\Q$ admits an action of 
$(\C^*_\Q)^n=\C^*_\Q\times\cdots\times\C^*_\Q$ ($n$ factors).
This action has a dense principal orbit. 
\end{proposition}
The action of $(\C^*_\Q)^n$ on $X_\Q$ is equivariant in the following sense: for each $\Lambda$ subgroup of finite index
of $\Z^n\cong\pi_1(U)$ there is a map $p_\lambda:X_Q\to{X_\Lambda}$ which is equivariant with respect to the actions of
$(\C^*_\Q)^n$ and  $(\C^*)^n$ on $X_\Q$ and $X$, respectively. Therefore the proalgebraic toric completion $X_\Q$ is an equivariant compactification of $(\C^*_\Q)^n$.

\section{Proalgebraic completion of a toric variety}\label{proalgebraic_completion}

For preliminaries on toric varieties we refer to the classical references \cite{Fulton}, \cite{Cox}. For an emphasis on the functorial properties of the theory we refer to \cite{KLMV}. For proalgebraic completions we refer the reader to \cite{SGA1}.

Given a cone $\sigma\subset\R^{n}$, consider the dual semigroup $S^\sigma$ of $\sigma\cap \Z^{n}$. Denote by $\mbox{Hom}_{sg}$ the space of semigroup morphisms. By Gordon's Lemma, $S^\sigma$ is a finitely generated semigroup. 

\noi Applying the covariant functor $\mbox{Hom}_{sg}(S^\sigma, \ \ )$ on the inverse system in diagram \eqref{diagram_III} we get the following inverse system:
\begin{equation}\label{diagram_IV}
\xymatrix{
\ldots\ar[r] & X_\sigma-\{p\} \ar[rr]^{(p_{m,n})_{*}}	&	& 	X_\sigma-\{p\} \ar[r] & \ldots\\
\ldots\ar[r] & (\Z/n\Z)^{r} \ar@{^{(}->}[u] \ar[rr]^{(p_{m,n})^{r}} &	&	(\Z/m\Z)^{r} \ar@{^{(}->}[u] \ar[r] & \ldots}
\end{equation}
whose inverse limit is the fiber of the following fibration:
\begin{equation}\label{diagram_V}
\xymatrix{\hat{\Z}^{r} \ar@{^{(}->}[r] & X^{\sigma}_\Q -\{cusp\} \ar[d]_{\pi_1}\\
& X^{\sigma}-\{p\} }
\end{equation}
where $X^{\sigma}:= \mbox{Hom}_{sg}(S^\sigma, \C )$ is the affine toric variety, $p:= \mbox{Hom}_{sg}(S^\sigma, (0) )$ is the special point and $(p_{m,n})_{*}:= \mbox{Hom}_{sg}(S^\sigma, p_{m,n})$ is post-composition by $p_{m,n}$. Here $r$ is the rank of $S^\sigma$ as a semigroup and we have defined $$X^{\sigma}_\Q:= \mbox{Hom}_{sg}(S^\sigma, \C_Q )$$
a ramified covering over the special point $p:= \mbox{Hom}_{sg}(S^\sigma, (0) )$. This is the \textsf{proalgebraic completion} of the affine variety $X^{\sigma}$ with a cusp over $p$. In the same fashion as before, the $\mbox{Hom}$ functor on the complex version of \eqref{diagram_II} gives a densely immersed copy of $\C^{r}$ into $X^{\sigma}_\Q -\{cusp\}$. This will be called the baseleaf and translations of it give a lamination of $X^{\sigma}_\Q -\{cusp\}$.

Now we extend the previous affine case construction to a toric variety as follows: Given a fan $\mathcal{F}$, consider it as a direct system with respect to the canonical inclusions. Define the \textsf{toric functor}:
$$F^{\mathcal{F}}:= \lim_{\underset{\sigma\in\mathcal{F}}{\longrightarrow} }\ \mbox{Hom}_{sg}\left(S^\sigma,\ \ \right).$$
The toric functor on the complex plane semigroup gives the usual toric variety associated to the fan $\mathcal{F}$:
$$X^{\mathcal{F}}:= F^{\mathcal{F}}\left(\C\right).$$
Instead of regular complex points, now we take the adelic $\C_{\Q}$ ones and define:
$$X^{\mathcal{F}}_{\Q}:= F^{\mathcal{F}}\left(\C_{\Q}\right).$$
By definition of the toric functor, it is clear that it is covariant and commutes with the inverse limit. Then we have:
$$X^{\mathcal{F}}_{\Q}= F^{\mathcal{F}}\left( \lim_{\underset{p_{n,m}}{\longleftarrow} }\ \C   \right)\cong
\lim_{\underset{q_{n,m}}{\longleftarrow} }\ X^{\mathcal{F}}$$
where we have defined $q_{n,m}:= F^{\mathcal{F}}(p_{n,m})$. In the next section, an explicit description of the map $q_{n,m}$ is given and it will be clear that it is a finite covering ramified over the divisors specified by the rays of the fan. In particular, they are toric covers and $X^{\mathcal{F}}_{\Q}$ is the proalgebraic completion of $X^{\mathcal{F}}$.

Taking the first canonical projection, we have a map:
$$X^{\mathcal{F}}_\Q\overset{q_1}{\longrightarrow} X^{\mathcal{F}}$$
ramified over the divisors specified by the rays of the fan.
\begin{definition} \label{solenoidization}
The map $q_1$ is called the \textsf{solenoidization} of  $X^{\mathcal{F}}$.
\end{definition}

In the case of a projective toric variety, the fiber is better pictured in the Delzant polytope associated to the fan. This can be described as follows: There is an equivariant action respect to the solenoidization of the adelic solenoid and the circle on $\C_\Q$ and $\C$ respectively such that:
\begin{equation}\label{quotient}
\C_\Q\rightarrow\C_\Q/S^{1}_\Q\cong \C/S^{1}\cong [0,+\infty)
\end{equation}
where the actions are free and properly discontinuous except at the origin. In particular, the fiber of the map \eqref{quotient} is a point at the origin and the adelic solenoid $S^{1}_\Q$ at the other points. By functoriality, this action extends to the proalgebraic completion of the toric variety and the orbit space is homeomorphic to the \textsf{Delzant polytope} of the fan:
$$\mbox{Delzant}(\mathcal{F}):= F^{\mathcal{F}}\left([0,+\infty)\right).$$

Applying the toric functor to the diagram \eqref{quotient} we get the \textsf{moment map}:
$$\mbox{M}:X^{\mathcal{F}}_\Q\rightarrow \mbox{Delzant}(\mathcal{F}).$$
Because of \eqref{diagram_V}, we immediately have:

\begin{proposition}\label{Delzant_prop}
Consider a projective toric variety and its proalgebraic completion. The fiber of the moment map of the completion over the points lying on an $m$-dimensional face of the Delzant polytope is isomorphic to the $m$-dimensional solenoidal torus $\left(S^{1}_\Q\right)^{m}$. In particular, the proper codimension one faces of the polytope correspond to the singular divisors and there are as many cusps as vertices of the polytope.
\end{proposition}

As an example, consider the projective space $\C P^{2}$ and its proalgebraic completion $\C P^{2}_\Q$. The respective Delzant polytope is pictured in Figure \ref{Fiber_Dimension}. The vertices correspond to cusp singularities, at the edges the moment map's fiber is $S^{1}_\Q$ and in the interior the fiber is $\left(S^{1}_\Q\right)^{2}$.

\begin{figure}
\begin{center}
  \includegraphics[width=0.35\textwidth]{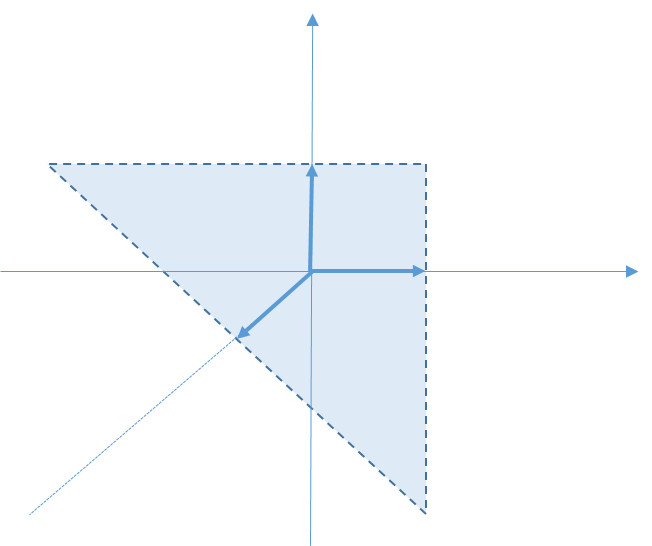}
  \end{center}
  \caption{\Small Delzant polytope of $\C P^{2}_\Q$. The vertices correspond to cusp singularities, at the edges the moment map's fiber is $S^{1}_\Q$ and in the interior the fiber is $\left(S^{1}_\Q\right)^{2}$.}\label{Fiber_Dimension}
\end{figure}

\subsection{Proof of theorem \ref{main_vector}: $\protect\mbox{Vect}_r\protect\left(X^{\mathcal{F}}_\Q\protect\right)\protect\cong \protect\varinjlim_{i\in I} \protect\mbox{Vect}_r\left(X^{\mathcal{F}}_i\right)$}\label{vector_section}

Theorem \ref{main_vector} follows immediately from the following Lemma.

\begin{lemma}\label{main_Lemma}
Consider a compact toric variety $X$ and the canonical maps $\pi_n: X_\Q\rightarrow X_n$ for every index $n$ in $I$. Consider a rank $r$ vector bundle $E\rightarrow X_\Q$. Then, there is an index $l$ and a rank $r$ vector bundle $E_l\rightarrow X_l$ such that its pullback by $\pi_l$ is isomorphic to $E$. \end{lemma}
\dem
Recall that $I$ is the index set of the projective system of finite coverings $X_i \rightarrow X$ branched over the invariant divisor under the torus action.

Consider the classifying map $e:X_\Q\rightarrow Gr(m,r)$ of the vector bundle $E$ with $m$ some natural number; i.e. The vector bundle $E$ is the pullback of the universal bundle by $e$. Consider a Whitney embedding $W:Gr(m,r)\hookrightarrow \R^{N}$ with $N$ some natural number and consider a, $\varepsilon$--tubular neighborhood $p:M\rightarrow Gr(m,r)$ with $M\subset \R^{N}$.

For every index $n$ in $I$ and every $x\in X_n$ we have the fiber $F_{n,x}$ of the bundle $\pi_n:X_\Q\rightarrow X_n$. This fiber is isomorphic (non canonically) to the adelic completion of the integers $\hat{\Z}$, a Hausdorff compact abelian topological group. In particular there is a unique normalized Haar measure $\mu_{n,x}$ on every fiber $F_{n,x}$.

Define the map $\zeta:= W\circ e:X_\Q\rightarrow \R^{N}$ and $\phi_{n,x}:F_{n,x}\rightarrow\R^{N}$ as the pre-composition of $\zeta$ with the fiber. We summarize our definitions in the  following commutative diagram:

$$\xymatrix{ & E \ar[r] \ar[d] \ar@{}[dr]|{\mbox{{\bf p.b.}}} & \mathcal{G}(m,r)\ar[d] & \\
\phi_{n,x}:\ F_{n,x} \ar@{^{(}->}[r] & X_\Q \ar@{{}{ }{}}@/_1pc/[rr]|{\mbox{\#}}\ar[d]_{\pi_n}\ar[r]^{e} \ar@/_2pc/[rr]_\zeta & Gr(m,r) \ar@{^{(}->}[r]^{\ \ \ W} & \R^{N} \\
 & X_n & & }$$

For every index $n$ in $I$, define the average map $f_n:X_n\rightarrow\R^{N}$:
$$f_n(x):=\int_{F_{n,x}}\phi_{n,x}\ d\mu_{n,x}.$$
By continuity, it is clear that $f_n\circ \pi_n$ tends to $\zeta$ pointwise and since $X_n$ is compact, there is a natural $l$ such that the image of $f_l$ is contained in the $\varepsilon$--tubular neighborhood $M$. Define the homotopy $H:X_\Q\times [0,1]\rightarrow Gr(m,r)$:
$$H(y,t):=p\left(t\zeta(y)+(1-t)f_l(\pi_l(y))\right).$$
Define the classifying map $e_l:= p\circ f_l$ and consider its corresponding rank $r$ vector bundle $E_l\rightarrow X_l$. Then, the classifying map $e$ is homotopic to $e_l\circ \pi_l$ by $H$; i.e. the following diagram $2$--commutes in the $2$--topological category with homotopies:

$$\xymatrix{X_\Q \ar@{{}{ }{}}@/_1pc/[rr]|{\mbox{{\it 2}\#}}\ar[d]_{\pi_l}\ar[rr]^{e} & & Gr(m,r) \\
  X_l \ar@/_1pc/[urr]_{e_l}& & }$$
and this finishes the proof.
\fdem

\section{Holomorphic vector bundles over $\PQ$.}\label{hol-vect-bun}

\subsection{Laurent-Puiseux series in $\C^*_\Q$}\label{lp} 
The abelian group $\CQ$ is isomorphic to the direct product of $S^{1}_{\Q}$   
with the multiplicative group of the positive reals i.e. 
$\CQ=\SQ\times{\R}^*$. This follows because it is an inverse limit of groups over $\C^*$ and $\C^*\simeq\bs^1\times\R^*$. If $z\in\CQ$ using this factorization 
$z=(\mathbf{u},t)$, with $\mathbf{u}\in\SQ,\, t\in\R^*$, we will write $z=t\mathbf{u}$ (i.e. we introduce ``\emph{polar coordinates}'' in $\CQ$). 

\begin{definition}\label{disks}
If $z\in\CQ$ we define the absolute value $|z|$ of $z$ as follows: $|z|=t$ if $z=t\bu$ (in polar coordinates). If one considers 
$\C_\Q=\CQ\cup{\boldsymbol{0}}$ (the inverse limit of
the branched self-coverings $z\mapsto{z^n}$ of $\C$ with canonical projection 
$p_1:\C_Q\to\C$ and $\boldsymbol{0}=p_{1}^{-1}\left\{0\right\}$ (see definition \ref{C*QandCQ}) then we define 
$|\boldsymbol{0}|=0$.

Let $r>0$. Define the open and closed disks of radius $r$ centered at the origin as follows:
\[
D^+_\Q(r)=\left\{q\in\C_Q\,:\, |q|<r \,\text{or} \,\, q=\boldsymbol{0} \right\},\quad \bar{D}^+_\Q(r)=
\left\{q\in\C_Q\,:\, |q|\leq{r}\, \text{or} \,\, q=\boldsymbol{0} \right\}
\]
Analogously, one defines the open and closed disks of radius $r$ centered at $\boldsymbol{\infty}\in\PQ$, where, 
$\boldsymbol{\infty}=p_1^{-1}(\left\{\infty\right\})$ and $p_1:\PQ\to\C P^1$ is the canonical projection, as follows:
\[
D^-_\Q(r)=\left\{q\in\C_Q\,:\, |q|>r \,\,\text{or}\,\,q=\boldsymbol{\infty} \right\},\quad \bar{D}^-_\Q(r)=
\left\{q\in\C_Q\,:\, |q|\geq{r} \,\,\text{or}\,\,q=\boldsymbol{\infty} \right\}.
\]
\end{definition}
 One has the decomposition of the adelic projective line into two closed ``hemispheres":
\[
\PQ=\bar{D}^+_\Q(1)\cup\bar{D}^-_\Q(1), \quad 
\text{with ``equator''}\,\, \bar{D}^+_\Q(1)\cap\bar{D}^-_\Q(1)=\SQ. 
\]
 \noi The south and north poles are $\boldsymbol0$ and $\boldsymbol\infty$, respectively.

 We recall that the Pontryagin dual $K$ of the compact abelian group 
 $S^{1}_{\Q}$ is the additive group of rational numbers $(\Q,+)$:

\[
K\overset{def}= 
\left\{f:S^{1}_{\Q}\to\bs^1\, : \,\text{f is a continuous homomorphism}  \right\}\simeq\Q,
\]

The group of additive endomorphisms of $(\Q,+)$ is $\Q$ acting by multiplication.
Hence the group of isomorphisms of $(\Q,+)$ is $\Q^*=\Q-\left\{0\right\}$ and 0
corresponds to the trivial endomorphism. Hence by Pontryagin duality the group of automorphisms
of $\SQ$ is isomorphic to $\Q^*$.
For each rational number $q\in\Q^*$, using the dual of the automorphism of
$\Q$, we obtain an automorphism 
$\xi_q:S^{1}_{\Q}\to\SQ$, corresponding to $q$ and, in polar coordinates, we define the map
$\hat\xi_q:\C^*_\Q\to\C^*_\Q$ by the formula:
$\hat\xi_q=\hat\xi_q(\mathbf{u},t)= t^q\xi_q(\mathbf{u})$ if $z=t\mathbf{u}$. To simplify the notation we will write briefly $\hat\xi_q(z)=z^q$. 

Therefore, ``taking a rational power''
of an element in $\CQ$ makes sense. In addition, and unlike the case of taking integral powers of $\C$, 
we have: 
\begin{enumerate}
\item For any rational $q\neq0$, the map $z\mapsto z^q$ 
is a continuous isomorphism of $\CQ$.
\medskip
\item For $q\neq0$ the inverse of the map $z\mapsto{z^q}$ is the map 
$z\mapsto{z}^\frac1q$ .
\end{enumerate} 
\begin{definition}\label{chc*} For each $q\in\Q$, define the homomorphism  
$\hat\chi_q:\CQ\to\C^*$ given by the formula
 $\hat\chi_q(z)=t^q\chi_q(\bu)$,
where $\chi_q$ is the character of $\CQ$ corresponding to $q$  and $z=t\bu$
in polar coordinates.
A series $S$ of the form $S=\underset{q\in\Q}\sum\,a_q$ where $a_q\in\C^n$ means
$\overset{\infty}{\underset{n=1}\sum}\,a_{q_n}$ for some numbering of the rationals 
$\Q=\left\{q_1,q_2,\cdots,\right\}$. If the series is absolutely convergent
i.e., the series $\overset{\infty}{\underset{n=1}\sum}\,|a_{q_n}|$ converges,  then the series $S$ converges to the same vector independently of the numbering. We say that $S$ converges absolutely.
\end{definition}
\begin{definition}[{\bf Holomorphic functions}]\label{holo} A continuous function $f:D^+_\Q(r)\rightarrow\C^n$ is said to be holomorphic
if it can be developed as a power series which is absolutely convergent and it is locally uniformly convergent (i.e. convergent on on compact sets):
\[
f(z)=\sum_{q\in\Q,\,q\geq0}\,a_q\hat\chi_q(z)=\sum_{q\in\Q,\, q\geq0}\,a_qt^q\chi_q(\bu),\quad f(\boldsymbol{0})=a_0
\]
where $a_q\in\C^n$, $\chi_q$ is the character corresponding to 
$q\geq0$ and $z=t\bu$, $0\leq{t}<{r}$. 

\noi Analogously, A continuous function $f:D^-_\Q(r)\rightarrow\C^n$ is said to be holomorphic
in the disk centered at $\boldsymbol{\infty}$ if it can be developed as an absolutely and locally uniformly convergent power series:
\[
f(z)=\sum_{q\in\Q,\,q\leq0}\,a_q\hat\chi_q(z)=\sum_{q\in\Q,\, q\leq0}\,a_qt^q\chi_q(\bu),
\]
where $a_q\in\C^n$, $\chi_q$ is the character corresponding to 
$q\geq0$ and $z=t\bu$, $r<t<\infty$, where $f(\boldsymbol{\infty})=a_0$. 
\end{definition}

\begin{definition}For real numbers $r$, $R$ with $0\leq{r}<{R}\leq\infty$ define the solenoidal open annulus $A(r,R)$
as follows $A(r,R)=\left\{q\in\CQ\,:\, r<|q|<R \right\}$. Note that 
$A(0,\infty)=\C^*$
\end{definition}
\begin{definition}\label{holoannulus}
Let $f:A(r,R)\to\C^n$ then f is said to be holomorphic if f can be expressed
as an absolutely and locally uniform convergent series:
\[
f(z)=\sum_{q\in\Q}\,a_q\hat\chi_q(z)=\sum_{q\in\Q}\,a_qt^q\chi_q(\bu),
\]
where $a_q\in\C^n$, $\chi_q$ is the character corresponding to 
$q\geq0$ and $z=t\bu$, $r<t<R$.

\noi {\it The series is required to be convergent only for points in the annulus}. 
Such a series will be called a {\bf Laurent-Puiseux series} (since $q<0$ is allowed). 

\end{definition}
\begin{definition}\label{LPPoly}
If the sum is of the form
\[
f(z)=\sum_{q\in{F}}\,a_q\hat\chi_q(z)=\sum_{q\in{F}}\,a_qt^q\chi_q(\bu),
\]
where $F\subset\Q$ is a finite set, then $f$ is called a {\bf finite Laurent-Puiseux series} if the sum is over a finite set of nonnegative rational is called
{\bf Laurent-Puiseux polynomial}.
 
\end{definition}

\subsection{Winding numbers}\label{wn} Let $0\leq{r}<1<R\leq\infty$. Let $f:A(r,R)\to\C^*$ be a holomorphic map then, since $A(r,R)$ retracts strongly onto $\SQ$ (as is easily seen using polar coordinates) and $\C^*$ retracts strongly onto the circle $\bS^1$ we see that
the homotopy class of $f$ is determined by the homotopy class of the map $\hat{f}$
 from $\SQ$ to $\bS^1$ given by the formula $\hat{f}(\bu)=\frac{f(\bu)}{|f(\bu)|}$.
 By Scheffer's theorem \cite{Sc} $\hat{f}$ is homotopic to a character  
 $\chi_q:\SQ\to\bS^1$ corresponding to a rational number $q\in\Q$.
 \begin{definition}\label{hwn} The rational number $q$ in the previous paragraph is called
 the {\bf winding number} of $f$ and it is denoted $w(f)$.
 \end{definition}
 \begin{definition}\label{winding_annulus}For any $r$, $R$ such that $0\leq{r}<{R}\leq\infty$ we can define the winding number of a continuous map 
 $f:A(r,R)\to\C^*$. 
Consider the restriction  of $f$ to the set 
$r'\SQ=\left\{r'\bu\,:\, \bu\in\SQ \right\}$ for any $r'$ with $r<r'<R$.  
This set is homeomorphic to $\SQ$ and therefore the homotopy class of this restriction determines
a character $\chi_q$ (independent of $r'$). Then $q\in\Q$ is called the winding number of $f$ and is is also denoted $w(f)$.
\end{definition}

\begin{remark} \label{rmkwinding2}

 Let  $f:\SQ\to\C^*$
be a map with an absolutely convergent Fourier series:
\[
f(\bu)=\sum_{q\in\Q}\,a_q\chi_q(\bu),\quad \sum_{q\in\Q}\,|a_q|<\infty. 
\]
Since $f$ is continuous it has a winding number $q\in\Q$. Then
$f$ can be factored in a unique way as
$f(\bu)=\chi_q(\bu)g(\bu)$ where $g(\bu)\neq0$, $\forall\bu\in\SQ$ and
the winding number of $g$ is 0. Furthermore $g$ can be developed into
a Fourier series with absolutely convergent coefficients.

  Analogously, if $f:A(r,R)\to\C^*$ is a holomorphic function
in the annulus $A(r;R)$ with Laurent-Puiseux development 
\[
f(z)=\sum_{q\in\Q}\,a_q\hat\chi_q(z)=\sum_{q\in\Q}\,a_qt^q\chi_q(\bu), \quad z=t\bu, \,a_q\in\C,
\]
\noi then there exists a holomorphic function $g$ and a unique character $\chi_q$ such that
$f(z)= \hat\chi_q(z)g(z)$ where $g(z)\neq0\,\,\forall z\in{A(r,R)}$ and the winding number of $g$ vanishes. 

\noi In both cases $g$ is obviously defined as follows: $g(z)=\hat\chi_{_{-q}}(z) {f(z)}$. 
\end{remark}

\begin{proposition}\label{Liouville} {\bf [Liouville's theorem and maximum principle for $\PQ$]} We recall that $\C^*_\Q$ is a Riemann surface lamination.  Its leaves are densely embedded copies of $\C$. In addition, 
each of the punctured open disks $D^+_\Q(r)-\left\{\boldsymbol0\right\}$ and
$D^-_\Q(r)-\left\{\boldsymbol\infty\right\}$ is a lamination by complex surfaces
which are densely embedded copies of the hyperbolic plane. The restriction
of a holomorphic map on one of the disks is a holomorphic
map on each leaf. It follows that if $f:\PQ\to\C$ is holomorphic when restricted
to two disks $D^+(r)$ and $D^-(R)$ ($r>R$) then $f$ must be a constant since
it is holomorphic and bounded when restricted to a leaf of $\PQ$, which is a copy of $\C$, so by Liouville's theorem is constant on each leaf.  Since each leaf is dense and $f$ is continuous $f$ is constant. Similar arguments imply:
if $f:{D}^+_\Q(1)\to\C$ attains its maximum in a point $\bz_0\in{D}^+_\Q(1)$
then $f$ is constant.

\end{proposition}

\subsection{Definition of holomorphic vector bundle over the adelic 
projective line $\C P^{1}_\Q$}\label{vb} 
A theorem that will be used in this and the following sections is Scheffer's theorem \cite{Sc}:
\begin{theorem}[\bf Scheffer's Theorem]\label{scheffer}
Let $\mathfrak{G}$ be a compact connected topological group, and let 
$\mathfrak{H}$ be a locally compact abelian topological group. Then every 
continuous map $f:\mathfrak{G}\to\mathfrak{H}$ is 
homotopic to exactly one continuous homomorphism  $\hat{f}:\mathfrak{G}\to\mathfrak{H}$ and 
the homotopy can be chosen to preserve the identities.
\end{theorem}
 First let us consider the topological facts underlying Picard's theorem for complex vector bundles of complex dimension one.
We recall that $\PQ$ is topologically
homeomorphic to the suspension $\Sigma(\SQ)$ of $\SQ$.  By standard facts about
classification of vector bundles we know that each complex line bundle over $\PQ$
is topologically equivalent to the pullback of the tautological complex line bundle $V$ over $\C{P}^\infty$, $P:V\to \C{P}^\infty$. Therefore complex vector bundles over
$\PQ$ are in one-to one correspondence with the set (actually a group)
of homotopy classes $[\PQ, \C P^\infty]$ of maps from $\PQ$ to $\C P^\infty$.
Actually, since $\PQ$ is 2-dimensional every such bundle is obtained as
the pullback of the Hopf bundle: $h:\mathcal{O}_{\C P^1}(1)\to\C P^1$, under a continuous map
$P:\PQ\to\C P^1\simeq{\bS^2}$ so the set of equivalence classes of \emph{topological} complex line bundles is in one-to-one correspondence with homotopy classes $[\PQ,\C P^1]$, of maps from $\PQ$
to $\bS^2$. Any such map is the suspension 
$\Sigma(g): \Sigma{\SQ}\simeq\PQ \to\Sigma{\bS}^1\simeq\bS^2$ of a map $g:\SQ\to\bS^1$.

Therefore, by Scheffer's theorem \ref{scheffer}, $g$ is homotopic to
a character of $\SQ$. Therefore, 
$[\PQ,\C P^1]\simeq[\SQ,\bS^1]\simeq\Q^*$. 
From this we obtain:
\begin{theorem}{\bf Topological Picard Theorem}\label{Picard} The group, under tensor product, of equivalence classes of (topological) complex line bundles over $\PQ$ is isomorphic to $\Q$. The rational number $q$ associated to a complex line bundle $L_q$ is the first Chern number of the line bundle. The first Chern number of the tensor product ${L_{q_1}}\otimes{L_{q_2}}$
is $q_1+q_2$.
\end{theorem}
 Let us define rigorously what we mean by a holomorphic vector bundle 
over $\PQ$. Given real numbers $r, R$ such that $0\leq{r}<1<{R}\leq\infty$ we define the open covering 
of $\PQ$ by the two open disks $D^+_\Q(R)$ and $D^-_\Q(r)$ centered at
$\boldsymbol0$ and $\boldsymbol\infty$, respectively. The intersection
$D^+_\Q(R)\cap{D^-_\Q(r)}$ is the open annulus $A(r,R)$. 

\noi Then
{\bf by definition} any $n$-dimensional holomorphic vector bundle
is given by the single holomorphic cocycle:
\[
f_{_{(r,R)}}(z): A(r,R)\to{\text{GL}}(n,\C).
\]
\medskip
The cocycle corresponds to the open covering of $\PQ$ by the two given open disks. 

\noi We have required $f_{(r,R)}$ to be a holomorphic map from the annulus
$A(r,R)$ to ${\text{GL}}(n,\C)$ which is an open set in the complex vector space
on $n\times{n}$ complex matrices $M_n(\C)$. Then, by definition \ref{holoannulus} $f_{(r,R)}$ can be developed
as a Laurent-Puiseux matrix-valued series
\[
f_{(r,R)}(z)=\sum_{q\in\Q}\,a_q\hat\chi_q(z)=\sum_{q\in\Q}\,a_qt^q\chi_q(\bu), \quad r<t<R,\, \bu\in\SQ,\quad a_q\in M_n(\C).
\]
 
 \begin{definition}\label{hvb} The 
  holomorphic vector bundle $\pi:E\to\PQ$  with fiber $\C^n$
obtained from the disjoint union of $D^+_\Q(R)\times\C^n$ 
and $D^-_\Q(r)\times\C^n$ by identifying $(z,W)$
with $(z,f_{(r,R)}(z)(W))$ for $z\in{A(r,R)}$ is called the vector bundle corresponding to the cocycle $f_{(r,R)}$ over each disk.
\end{definition}
 Thus the construction of the bundle given in definition \ref{hvb} 
 is the standard one using only one clutching function. The fibration $\pi$ in definition \ref{hvb}  
is holomorphically trivial over each of the disks $D^+_\Q(R)$ and $D^-_\Q(r)$ 
 i.e. there are $n$ linearly independent holomorphic sections 
 (holomorphic in the sense of definition \ref{holo}). 

The holomorphic equivalence of two cocycles $f_{(r,R)}$
 and $f_{(r',R')}$ is the usual one. The cocycle determined by $f_{(r,R)}$
is equivalent to the cocycle determined by $f_{(r',R')}$ if there
exists two holomorphic maps $f_1, f_2:A(r'', R'')\to{\text{GL}}(n,\C)$ where $r''=\max\left\{r,r'\right\}$, $R''=\min\left\{R,R'\right\}$
such that $f_1f_{(r,R)}f_2^{-1}=f_{(r',R')}$. Equivalent cocycles define equivalent bundles.

\medskip
Theorem \ref{Picard}, in the topological category, is also a consequence of the functoriality implicit in Theorem \ref{main_vector}.
Consider the complex projective line $\C P^1$ and its 
proalgebraic completion $\PQ$.  The index set is the set of natural numbers $I=\N$. 
Indeed, it is clear that
\begin{equation}\label{identity}
\I_{\C P^1}\left(\frac{n'}{n}m\right)\,\cong\, p_{n, n'}^*\,\I_{\C P^1}(m)
\end{equation}
for every pair of naturals such that $n|n'$ and every integer $m$. 

We clearly have the following identity:
\begin{equation}\label{pb}
\I_{\C P^1_\Q}\left(\frac{m}{n}\right)= \pi_n^{*}\,\I_{\C P^1}(m)
\end{equation}
where $\pi_n:\C P^1_\Q\rightarrow\C P^1$ is the canonical projection of the inverse limit to level $n$ of the projective system.

Because of the additivity of the characters, identity \ref{identity} and the fact that the pullbacks of the projections are equivariant with the tensor product up to isomorphism, we have
$$\I_{\C P^1_\Q}(q)\otimes \I_{\C P^1_\Q}(q')\cong \I_{\C P^1_\Q}(q+q')$$
for every pair of rational numbers $q$ and $q'$. In particular, \ref{pb} is well defined. This is obviously a holomorphic line bundle according to our
definition \ref{hvb}. 
\begin{remark} Here we considered line bundles which are pullback of the Hopf
bundle $\mathcal{O}(1)$. These are holomorphic in the sense defined in \ref{hvb}
above. 

\end{remark}

\begin{proof}[Proof of Theorem \ref{Picard} using Theorem \ref{main_vector}] 
We saw that every topological complex line bundle on $\C P^1_\Q$ ($\C P^1$) is in one to one correspondence with a homotopy class in $[S^{1}_\Q, S^{1}]$ ($[S^{1}, S^{1}]$). Thus,  every line bundle is homotopic to a holomorphic line bundle. Then, by Theorem \ref{main_vector} for holomorphic line bundle categories
$$Pic\left(\C P^1_\Q\right)\cong \lim_{\underset{n\in\N}{\longrightarrow} }\,Pic\left(\C P^1_n\right)
\cong \lim_{\underset{n\in\N}{\longrightarrow} }\,\frac{1}{n}\Z\cong \Q$$
and we have the result.
\end{proof}

\section{Adelic loop groups}\label{adelic-loop-groups} 

In analogy with classical loop groups
(see, for instance \cite{PS}, \cite{MG, MG1} ), and following the dictum of Dennis Sullivan that $\SQ$ is a ``diffuse'' version of $\bS^1$, we have the following:
\begin{definition}\label{adelic_loop} Let $G$ be a connected Lie group.  
A continuous {\bf adelic loop} in $G$ is a continuous map $f:\SQ\to{G}$. The
{\bf adelic loop group} of $G$, denoted $\boldsymbol\Lambda_\Q(G)$, is defined as follows:
\[
\boldsymbol\Lambda^0_\Q(G)=\left\{f:\SQ\to{G}\,:\, f\,\text{is a continuous loop} \right\}.
\]
Then $\boldsymbol\Lambda^0_\Q(G)$ is a topological group under point-wise
multiplication and the compact-open topology.

\end{definition}

\begin{remark}\label{differentiable_loop} In definition \ref{adelicsmoothloopgroup} we will define the
group of ``differentiable'' loops ${\boldsymbol\Lambda}^{\infty}_\Q(G)$ consisting
of ``differentiable'' loops $f:\SQ\to{G}$ such that $f$ restricted to each leaf of the solenoidal
manifold $\SQ$ (a one dimensional copy
of $\R$ densely embedded in $\SQ$) is differentiable.
Differentiable based solenoidal groups loops form an infinite dimensional
Fréchet manifold (see, for instance  \cite{Fr},\cite{MG,MG1}, \cite{PS}). 
 However, in the following sections we will use loops with matrix coefficients
 in the Wiener algebra and $H^1$ Sobolev space (which is also a Banach algebra).
\end{remark}

We are mostly interested in the {\bf adelic loop group} $\boldsymbol\Lambda^0_\Q(\gln)$. Since
the inclusion $\text{U}(1)\hookrightarrow\gln$ induces an isomorphism of fundamental groups we have $\pi_1(\gln)=\Z$. Here $\text{U}(1)$ is included as the center of 
$\gln$
of scalar matrices so that $\gln/\text{U}(1) =P \gln$, and one has the fibration 
$p_1:\gln\to P\gln$ with fibre $\text{U}(1)$.
 
\begin{proposition} Any loop $f:\SQ\to\gln$ can be deformed to a map 
$\hat{f}:\SQ\to{\text{U}(1)}\simeq\bS^1$ into a fibre
of $p_1$. 
 \end{proposition}
 
 \begin{proof} The proof follows from the  fact that $\SQ$ has topological dimension one. In fact, one has the fibration $p_2:P\gln\to \C P^{n-1}$ and the  
loop $f$ can be deformed to a map that 
avoids the union $H$ of the pre-images of the coordinate hyperplanes under the projection $p_2\circ{p_1}$ (since $\SQ$ is of real dimension one). The complement 
$\gln - H$ is diffeomorphic to 
$\C^{n-1}\times{\text{U}(1)}$ where the $\text{U}(1)$-factors are fibers of $p_1$.  By theorem \ref{scheffer}, $\hat{f}$ is homotopic to a continuous group homomorphism to $\text{U}(1)\simeq\bS^1$ i.e,. a character. 
\end{proof}

\noi Therefore from this proposition we have:
\begin{proposition}\label{connectedcomponents}
 The group of homotopy
classes of loops in $\La^0_\SQ(\gln)$ is isomorphic to $\Q$. Thus the connected components of $\boldsymbol\Lambda^0_\Q(\gln)$ are in one-to-one 
correspondence to the additive group of $\Q$. The connected component of the identity corresponds to q=0.
\end{proposition}

\medskip
\begin{definition}[{\bf Determinant bundle of complex vector bundles over $\PQ$}]\label{determinant-bundle} 

\noi Let $\pi:V\to\PQ$ be a holomorphic vector bundle of rank $r\geq1$.  Suppose that
$V$ corresponds to the cocycle $f:A(r,R)\to\gln$. Let $d_f:A(r,R)\to\C^*$
be defined by composing with the determinant as follows $ d_f(z)=\text{det}\,f(z)$.
The loop $d_f$ is called the {\bf determinant loop}; $d_f$ determines a holomorphic cocycle and thus a holomorphic line bundle
denoted $\text{det}(V)$ on $\PQ$. 

\noi The bundle $\text{det}(V)$ is called the
{\bf determinant line bundle of $V$}. Its first Chern class is the rational 
number $w(d_f)$ which will be denoted $w(V)$. 
It is the winding number according to definition \ref{winding_annulus}. The connected component of the identity of the loop group $\boldsymbol\Lambda^0_\Q(\gln)$ consists of those 
vector bundles $V$ of rank $n$ such that $w(V)=0$. Of course $w$ is a topological invariant of these vector bundles. The determinant is multiplicative: 
$d_{f\cdot{g}}=d_f{d_g}$.
\end{definition}

\medskip

\begin{remark}[\bf Dictum]\label{dictum}Almost everything that holds for the standard loop group $\boldsymbol\Lambda^0(G)$, for loops with any regularity assumptions (differentiable, in Wiener apace, H\"older, Sobolev etc.), holds for the solenoidal group loop 
$\boldsymbol\Lambda^0_\Q(G)$, where $G$ is a connected Lie group.
\end{remark}

The {\it dictum} is due to the fact that the theory for standard loops, for instance in \cite{CG, Fr, MG, MG1, PS, Pr, SW}, use only the fact that $S^1$ is a compact. connected, abelian Lie group. In our case $\SQ$ 
is a {\bf connected pro-Lie group} (\cite{HM}), however $\SQ$ is not locally connected. Therefore, anything that is done in $\boldsymbol\Lambda^0(G)$ that 
does not use local connectedness can be applied to $\boldsymbol\Lambda^0_\Q(G)$.
One difference with the classical loop groups is that the Character group of $\SQ$ is the additive rationals $(\Q,+)$ instead of $\Z$, so in many of the statements 
$\Q$ takes the place of $\Z$. This enriches the theory.
In fact, to show why the dictum \ref{dictum} holds, we will follow many of the steps in Chapters 6, 7
and 8 in the book by Pressley and Siegel \cite{PS}.

\subsection{Loops in the Wiener Banach algebra}
\label{Wienerloops}
Let $\mathfrak{W}_\Q(n)$ denote the algebra of functions 
$f:S^{1}_{\Q}\rightarrow M_n(\C)$, for some fixed integer $n\geq1$,
such that $f$ can be expressed as a Laurent-Puiseux series of the form
\[ 
f(\bu)=\sum_{q\in\Q}\, A_q\,\chi_q(\bu),  
\quad \bu\in\SQ,\,\,A_q\in M_n(\C)\,\, \text{the ring of}\,\, n\times{n}\, 
\text{complex matrices}.
\]
and
\[
\sum_{q\in\Q}\,|A_q|<\infty
\]
Then $\mathfrak{W}_\Q(n)$ is a Banach algebra with norm 
\[
|f|=\sum_{q\in\Q}\,|A_q|.
\]
Multiplication is the standard multiplication of matrix functions. 
This Banach algebra it is called the {\bf Wiener Algebra}.

Gelfand's proof (\cite{Ge1, Ge2}) of Wiener-Lévy's  theorem (\cite{K, Le, W}) is valid for any Banach algebra (\cite{BP}) and therefore applies to the group of units 
$\boldsymbol\Lambda_W(n)$ of the Banach algebra
$\mathfrak{W}_\Q(n)$,
and it implies that $f$ is a unit if and only if the image of $f$ lies in
$\text{GL}(n,\C)$ i.e. $f(\bu)$ is an invertible matrix for all 
$\bu\in S^{1}_{\Q}$.
The set
\[
\boldsymbol\Lambda_W(n)\overset{\small{def}}=
\left\{ f\in\mathfrak{W}_\Q(n)\,:\, f(\bu)\,\, 
\text{is invertible for all}\,\,   \bu\in S^{1}_{\Q}  \right\}
\]
is open in $\mathfrak{W}_\Q(n)$ and thus it is a 
{\bf Banach Lie group} (\cite{Mi}).  

\begin{definition}[Wiener loop]\label{Wl}
 An element $f\in\boldsymbol\Lambda_W(n)$ is called a {\bf Wiener loop}
 (or $l^1$-loop)
of the group $\boldsymbol\Lambda_W(n)$ (\cite{BD, CG, GK,PS}).
\end{definition}

Let 
$$\mathfrak{W}^+_\Q(n)=\left\{ f\in\mathfrak{W}_\Q(n) \,: \, 
f(\bu)=\sum_{q\in\Q}\, A_q\,\chi_q(\bu),\, u\in\SQ\,\,
\text{such that}\, A_q=0,\, \text{for}\,q<0 \right\},$$
and
$$\mathfrak{W}^-_\Q(n)=\left\{ f\in\mathfrak{W}_\Q(n) \,: \, 
f(\bu)=\sum_{q\in\Q}\, A_q\,\chi_q(\bu),\, \, u\in\SQ\,\,
\text{such that}\, A_q=0,\, \text{for}\,q>0 \right\},$$

$\mathfrak{W}^+_\Q(n)$ and $\mathfrak{W}^-_\Q(n)$ are closed subalgebras
of $\mathfrak{W}_\Q(n)$ and there is a splitting of closed subalgebras:

\[
\mathfrak{W}^+_\Q(n) + \mathfrak{W}^-_\Q(n)=\mathfrak{W}_\Q(n), \quad
\mathfrak{W}^+_\Q(n) \cap \mathfrak{W}^-_\Q(n)\simeq{M_n(\C)}
\]
\noi A function is in $\mathfrak{W}^+_\Q(n)$ if it can be extended to a 
continuous function which is holomorphic in the open disk $D^+_\Q(1)$:
if $f\in\mathfrak{W}^+_\Q(n)$ is given by the series
$$
f(\bu)=\sum_{q\in\Q}\, A_q\,\chi_q(\bu),\quad u\in\SQ,\,\text{with}\,\, A_q=0,\, \text{for}\,q\leq0 
$$
\noi define its extension $f^+$ to $D^+_\Q(1)$ by the formula
$$
f^+(z)=\sum_{q\in\Q}\, A_q\,\hat\chi_q(z)=
\sum_{q\in\Q}\, A_q\,t^q\chi_q(\bu)
,\, \, z=t\bu,\, 0\leq{t}\leq1,
$$
\noi Analogously, a function is in $\mathfrak{W}^-_\Q(n)$ if it can be extended to a continuous function which is holomorphic in the open disk $D^-_\Q(1)$:
if $f\in\mathfrak{W}^-_\Q(n)$ is given by the series
$$
f(\bu)=\sum_{q\in\Q}\, A_q\,\chi_q(\bu),\quad u\in\SQ,\,\text{with}\,\, A_q=0,\, \text{for}\,q\geq0 ,
$$
\noi define its extension $f^-$ to $D^-_\Q(1)$ by the formula
$$
f^-(z)=\sum_{q\in\Q}\, A_q\,\hat\chi_q(z)=
\sum_{q\in\Q}\, A_q\,t^q\chi_q(\bu)
,\, \, z=t\bu,\, 0\leq{t}\leq1,
$$
Hence, $\mathfrak{W}^+_\Q(n)$ 
is isomorphic to the algebra of continuous functions 
$f:\bar{D}^+_\Q(1)\to{M_n(\C)}$ such that $f$ has an absolutely convergent
Fourier series when restricted to $\SQ$ and it is holomorphic in 
the open disk $D^+_\Q(1)$. Also $\mathfrak{W}^-_\Q(n)$ 
is isomorphic to the algebra of continuous functions 
$f:\bar{D}^-_\Q(1)\to{M_n(\C)}$ such that $f$ has an absolutely convergent
Fourier series when restricted to $\SQ$ nd it is holomorphic in 
the open disk $D^-_\Q(1)$

Let $\boldsymbol\Lambda_W^+(n)$ and $\boldsymbol\Lambda_W^-(n)$ denote the group of units
of $\mathfrak{W}^+_\Q(n)$ and $\mathfrak{W}^+_\Q(n)$, respectively. In view
of the identification above, 
and Wiener-Lévy's theorem both  $\boldsymbol\Lambda_W^+(n)$ and 
$\boldsymbol\Lambda_W^-(n)$ 
consists of loops $f$ such that, when extended to the corresponding disks, their values consist of nonsingular matrices, since its elements are units. Therefore, both groups are Banach Lie groups {\cite{Mi}}.
\begin{proposition}\label{components-A_r} \mbox{}
\begin{enumerate}
\item The connected components, $\boldsymbol\Lambda_{W,q}(n)$ ($q\in\Q$) of 
$\boldsymbol\Lambda_W(n)$ are naturally indexed by the rationals, via
the winding number of the determinant map $\bu\mapsto\,\det f(\bu)$ for $f$ in a connected component. The connected component of the identity is 
$\boldsymbol\Lambda_{W,0}(n)$.
\item  $\boldsymbol\Lambda_W^+(n)\subset\boldsymbol\Lambda_{W,0}(n)$ and 
$\boldsymbol\Lambda_W^-(n)\subset\boldsymbol\Lambda_{W,0}(n)$, i.e., the winding numbers of the determinants of loops in each component are 0.
\end{enumerate}
\end{proposition}
\begin{proof}
(1) Indeed,
all the loops in a connected component are homotopic and therefore
the winding numbers of the associated determinant loops are constant.
Hence, we can index the connected component by the rational common winding number.
The connected component of the identity contains the trivial loop
that maps all $\SQ$ to the identity and therefore every loop in this component has vanishing winding number.

\noi (2) Since both $\bar{D}^+_\Q(1)$ and $\bar{D}^-_\Q(1)$ are contractible (they are cones over $\boldsymbol0$ and $\boldsymbol\infty$, respectively)
it follows that the loops $\gamma^+\in\boldsymbol\Lambda_W^+(n)$ and 
$\gamma^-\in\boldsymbol\Lambda_W^-(n)$ are homotopic to the constant matrices 
$\gamma^+(\boldsymbol0)$ and $\gamma^-(\boldsymbol\infty)$ in $\gln$, respectively. Thus
the winding numbers of the determinants of both of these loops is 0 so that they are both contained in the identity components $\boldsymbol\Lambda_{W,0}(n)$.
\end{proof}

%%%%
\section{Kähler structure and Energy Functional of adelic loop groups}\label{ksloops} 
\noi In this section we use the fact that $\SQ$ is a smooth 1-dimensional
lamination. Therefore we can use different notions of regularity, see for example
\cite{Mo, Su4, Ve}. \medskip
\begin{definition}\label{sobolev}
 Let $\bd\bu$ be normalized Haar measure on $\SQ$.
A function $f:\CQ\to\C^n\,$ is said to be
in the Sobolev space $H^1(\SQ,\C^n)$ if each component of $f$ is
a function $f_i\in H^1(\SQ)$, where $H^1(\SQ)$ is the Sobolev space  
of complex-valued functions defined as follows:
  \begin{align}
  H^1(\SQ)=\left\{g\in{L^2(\SQ,\bd\bu)}\boldsymbol{:}\, 
  f(\bu)=\sum_{q\in\Q}a_q\chi_q(\bu), 
  \quad\text{with}\,\, \sum_{q\in\Q} |a_q|^2(1+q^2)<\infty \right\}.
\end{align}
 \end{definition}

We recall that the exponential $\boldsymbol{Exp}$ in definition \ref{exp} is an epimorphism $\boldsymbol{Exp}:\R\times\hat{\Z}\to\SQ$ which is continuous.

 \begin{definition}\label{diffloops} If $M$ is a differentiable ($C^\infty$) manifold, 
then a map $h:\SQ\to{M}$ is said to be {\bf differentiable}
if the map $\tilde{h}=h\circ\boldsymbol{Exp}$ has the property
that for every $\boldsymbol{z}\in\hat{\Z}$
the restriction of $\tilde{h}$ to ${\R\times{\boldsymbol{z}}}$, defines
the map $t\overset{h_{\boldsymbol{z}}}\longmapsto\tilde{h}(t,\boldsymbol{z})$ from $\R$ to $M$ which
is a $C^\infty$ map (i.e. $\varphi\circ{h_{\boldsymbol{z}}}$ is smooth for any coordinate chart of $M$). 

\medskip 
\noi In the same way we can define functions of Hölder class $C^r\,\,, r>0 $. The derivative of $f':\SQ\to\R^n$ of class $C^1$ is defined as follows: $f'(t_0\bu_0)=\frac{d}{dt}h_{\bz_0}(t)|_{t=t_0}\,$ if 
$\boldsymbol{Exp}(t_{0},\bz_{0})=t_{0}\bu_0$.

\medskip
\noi Analogously, a function $h:\SQ\to{M}$ is said to belong to Sobolev space
$H^1(\SQ,M)$ if the map $\varphi\circ{h}$ belongs to Sobolev space 
$H^1(\SQ,\R^d)$ for every coordinate chart $\varphi:\mathcal{U}\to\R^d$ of $M$ (where $\dim_\R{M}=d$).

\end{definition} 

\noi The functions $f$ in $H^1(\SQ)$ are characterized by the properties:
  \begin{enumerate}
 \item $f\in{L_2(\SQ,\bd\bu)}$.
 \item $f'\in{L_2(\SQ,\bd\bu)}$ where $f'$ is the weak derivative
 of $f$ along the leaves of the lamination of $\SQ$. 

\end{enumerate}
 In \cite{TT} Proposition 1.1, is shown that if 
$f,g\in{W^{1,2}(\R)}$ then $||fg||_{_{W^{1,2}}} 
\lesssim ||f||_{_{W^{1,2}}}||g||_{_{W^{1,2}}}$. This implies, since 
$W^{1,2}(\bS^1)=H^1(\bS^1)$, that that the product of two functions
in $H^1(\bS^1)$ belongs to$H^1(\bS^1)$. The same proof applies to the
one-dimensional compact abelian group $\SQ$. One has:
 
\begin{proposition}\label{h1Banach-algebra} The product of two functions in $ H^1(\SQ)$ belongs
to $H^1(\SQ)$. Therefore $H^1(\SQ)$ is a Banach algebra with respect to the 
Sobolev norm.
\end{proposition}

\begin{definition}\label{adelicsmoothloopgroup}
Let $G$ be a compact Lie group with Lie algebra $\mathfrak{g}$. Let
$\Ld$
 be the group loop of smooth adelic loops on $G$:
 \[
 \Ld=\left\{\gamma:\SQ\to{G}: \gamma\,\, 
 \text{is smooth} \right\}. 
 \]
 The {\bf based adelic loop group} $\Lo$, consists 
 of loops $\gamma\in\Ld$ such that $\gamma(\boldsymbol1)=\boldsymbol{e}$
 where $\boldsymbol1$ is the identity element of $\SQ$ 
 (viewed as a multiplicative group) and $\boldsymbol{e}$ is the identity of
 $G$.  Both $\Ld$ and $\Lo$ are Fréchet manifolds modeled on the solenoidal loop space 
 of smooth maps from $\SQ$ to $\mathfrak{g}$. In a natural way $G$ acts on 
 $\Ld$ and therefore $\Lo$ can be identified with the homogeneous
 space $\Ld/G$, where $G$ is identified with the constant loops.
  \end{definition}
   We can also define 
  loop spaces modeled on Sobolev spaces. Following the steps of \cite{Fr, Pr, PS}
  we will consider loop spaces modeled on the Sobolev space $H^1(\SQ)$.

\begin{definition}\label{adelicsmoothloopgroup}
Let $G$ be a Lie group with Lie algebra $\mathfrak{g}$. Let
$\La^1(G)$
 be the group loop of  adelic loops on $G$ with values in $H^1(\SQ,G)$:
 \[
 \La^1(G)=\left\{\gamma:\SQ\to{G}: \, 
 \gamma\in{H^1(\SQ,G)} \right\}. 
 \]
 Let $\La^1_0(G)$ be the {\bf based adelic loop group}
 of loops $\gamma\in\La^1(G)$ such that $\gamma(\boldsymbol1)=\boldsymbol{e}$
 where $\boldsymbol1$ is the identity of $\SQ$ 
 (viewed as a multiplicative group) and $\boldsymbol{e}$ is the identity of
 $G$.  
  \end{definition}

  \medskip
  \noi {\bf In what follows we will mainly consider the case when $G=U(n)$}.
  Its Lie algebra $\mathfrak{U}(n)$, consist of $n\times n$ skew-Hermitian matrices which we will identify with $\R^{n^2}$. 
  
 \noi Consider the Sobolev based
 adelic loop space
  $\La_0^1(U(n))$. By proposition \ref{connectedcomponents} the group of connected components is isomorphic to $\Q$. 
  
  \medskip
  \noi {\bf For simplicity of notation we denote with the symbol $\bO$ the loop group $\La_0^1(U(n))$}.
  
  \medskip
  \noi Let 
  \begin{equation}
  \mathcal{L}\bO=\left\{f\in H^1(\SQ, \mathfrak{U})\,:\, f(\mathbf1)=0 
  \,\, =\text{the zero element of}\,\, \mathfrak{U}\right\}/\mathfrak{U}.
   \end{equation} \label{LOmega}
 \noi We call $\mathcal{L}\bO$ the based {\bf adelic loop algebra} of $\mathfrak{U}$.
  One has
 $$
 \mathcal{L}\bO={\mathfrak{u(n)}}\otimes H^1(\SQ)\simeq{H^1(\SQ,\R^{n^2})}.
 $$ 
 
 \noi The algebra
 $\mathcal{L}\bO$, with bracket 
 $$
 [g_{1}\otimes f_{1},g_{2}\otimes f_{2}]=[g_{1},g_{2}]\otimes f_{1}f_{2}
 $$
 
($g_{1}, g_{2}\in\mathfrak{u(n)},\, \, f_{1},f_{2}\in H^1(\SQ)$), is the tangent space at the identity of $\bO$ and the bracket is the Lie bracket of left
invariant vector fields on $\bO$. Hence $\mathcal{L}\bO$ is the the infinite dimensional Lie algebra of $\bO$ (see \cite{Ka} for details on infinite dimensional algebras).

 \begin{theorem}\label{kählerloop} (Compare \cite{Pr, Pr1, PS})   The based adelic loop space 
 $\bO$ is a connected, real analytic, Hilbert Lie group which in addition is an infinite dimensional, holomorphically homogeneous, complex Kähler manifold,
 modeled on the based loop algebra $\mathcal{L}\bO$.
 \end{theorem}
 
\begin{proof} 
By definition $\bO$ is connected. We refer to \cite{Pr1} Section 2 and
section 8.9 of \cite{PS} (see also \cite{Se}, \S 2,15, and Theorem 1.2 and also Appendix A in \cite{FU})
for the fact that $\bO$ is a real analytic Hilbert Lie group
modeled on the Hilbert space $\mathcal{L}\bO$. 

\noi First, let's remark that, in 
terms of Fourier series:
\begin{equation}\label{fourierH1}
\mathcal{L}\bO=\left\{f\in{H^1(\SQ,\mathfrak{U})}\, \boldsymbol{:}  
f(\bu)= \sum_{{q\in\Q,q<0}}\, a_q\chi_q(\bu)
+ \sum_{{q\in\Q,q>0}}\,a_q\chi_q(\bu)\,\boldsymbol{,}\,\,
a_q\in\mathfrak{U},\, a_q=\overline{a_{-q}}  \right\}.
\end{equation}

\noi The coefficient corresponding to $q=0$ vanishes since 
the loops are based (equation \ref{LOmega}).   

\medskip
\noi In what follows $\langle \,\cdot\,\,, \,\, \cdot\rangle_{_G}$ denotes a left and right-invariant Riemannian metric on $G$.

\medskip
\noi We can think of $\mathcal{L}\bO$ alternatively as the tangent space at
the identity loop $T_{\boldsymbol{e}}\bO$ or as left-invariant vector fields
on $\bO$.

\medskip
\noi The symplectic form, complex structure and Kähler metric are given explicitly as follows:

\begin{enumerate}
\item[\bf I.] {\bf Symplectic form.} Let
\begin{equation}
\omega(f,g)=\,\int_{\SQ}\langle f'(\bu),g(\bu)\rangle_{\mathfrak{U}}\, \bd\bu,
\end{equation}\label{sympform}

\noi where $\langle \cdot\,, \, \cdot \rangle_{\mathfrak{U}}$ is an
Ad-invariant inner product on $\mathfrak{U}$ so that $\omega$ is left-invariant. By equation \ref{LOmega}
the only $f\in{\mathcal{L}\bO}$ that has $f'\equiv0$ is the trivial 
loop that is identically 0. Then $\omega(f,f)>0$ if $f\neq0$, therefore
$\omega$ is non-degenerate. Integration by parts implies 
$\omega(f,g)=-\omega(g,f)$ so $\omega$ is skew.

 The proof that $w$ is closed is straight forward
and it is done in \cite{Pr1} pages 560-561. For completeness we include the adapted version of this argument here:

\begin{equation}
\begin{aligned} 
d\omega(X,Y,Z){} = & X\cdot\omega(Y,Z)- Y\cdot\omega(X, Z)+Z\cdot\omega(X, Y) \\
& -\omega([X,Y],Z)+ \omega([X,Z],Y)-\omega([Y, Z],X),
\end{aligned}
\end{equation}

where $X, Y,Z \in \mathcal{L}\bO$ are regarded as left invariant vector fields on $\bO$. As $\omega$ is left-invariant $\omega(Y, Z)$,  $\omega(X, Z)$ and $\omega(X, Y)$
are constant, hence each of the first three terms vanishes. Hence,

\begin{equation}
\begin{aligned} 
d\omega(X,Y,Z){} =  &\int_\SQ\, (-\langle [X',Y],Z\rangle_{\mathfrak{U}}
-\langle [X,Y'],Z\rangle_{\mathfrak{U}}+\langle [X',Z],Y\rangle_{\mathfrak{U}}
+\langle [X,Z'],Y\rangle_{\mathfrak{U}} \\
                     & -\langle [Y',Z],X\rangle_{\mathfrak{U}}
 -\langle [Y,Z'],X\rangle_{\mathfrak{U}} ) \, \bd\bu.
\end{aligned}
\end{equation}
Using the fact that $\langle \cdot\,,\, \cdot\rangle_{\mathfrak{U}}$ is invariant
under the adjoint function it follows that:
\begin{equation}
d\omega(X,Y,Z)=-2\left(\,\omega([X,Y],Z)+\omega(Z,[X,Y])\,\right)=0, \quad 
\forall\,X,Y,Z\in{\mathcal{L}\bO}.  
\end{equation}
\noi Hence $\omega$ is a left-invariant closed 2-form.

\medskip
\item[\bf II.] {\bf Complex structure.} As in \cite{Fr, Pr,Pr1, PS}, we define the almost-complex structure of the loop group $\bO$ using the characterization of its Lie algebra  $\mathcal{L}\bO $ in terms of 
Fourier series given in equation \ref{fourierH1}.
If
\[
f(\bu)= \sum_{{q\in\Q,q<0}}\, a_q\chi_q(\bu)
+ \sum_{{q\in\Q,q>0}}\,a_q\chi_q(\bu)\,\boldsymbol{,}\quad
a_q\in\mathfrak{U},\,\, a_q=\overline{a_{-q}},  
\]
\noi we define $J(f)$ as follows
\begin{equation}\label{complexstructure}
J(f)(\bu)=\sum_{q\in\Q,q<0}\mathbf{i}a_q\chi_q(\bu)
 -\sum_{q\in\Q,q>0}\mathbf{i}a_q\chi_q(\bu).
\end{equation}
\noi Obviously $J^2=-I$ so it is an almost-complex structure. As is pointed in \cite{Pr1}, the set of left-invariant vector fields of type $(0,1)$ is a subalgebra of $\mathcal{L} \bO$, and  the fact that both $\bO$  and $J$ are real analytic
implies by {\it Théorème 3.5} in \cite{Pe} that $\bO$ is a complex Hilbert manifold. The natural action on the left of $\bO$ on itself is by biholomorphisms, so $\bO$ is homogeneous.

\medskip
\item[\bf III.] {\bf Kähler metric}. If $f,h\in \mathcal{L}\bO$ have Fourier  series
\begin{equation}\label{f}
f(\bu)= \sum_{{q\in\Q,q<0}}\, a_q\chi_q(\bu)
+ \sum_{{q\in\Q,q>0}}\,a_q\chi_q(\bu)\quad (a_q=\overline{a_{-q}}),\, a_q\in\mathfrak{U}
\end{equation}
\noi and
\begin{equation}\label{g}
h(\bu)= \sum_{{q\in\Q,q<0}}\, b_q\chi_q(\bu)
+ \sum_{{q\in\Q,q>0}}\,b_q\chi_q(\bu) \quad (b_q=\overline{b_{-q}}), \,a_q\in\mathfrak{U}.
\end{equation}

The bilinear map $\boldsymbol{g}(f,h)=\omega(f,J(h))$ is definite positive. Indeed,
if $f$ and $g$ are given by the series \ref{f} and \ref{g}, respectively,
then:
\begin{equation}\label{kählermetric}
\boldsymbol{g}(f,g)=\sum_{q\in\Q}\,|q|\langle{a_q, b_q}\rangle_{_\mathfrak{U}}.
\end{equation}

\end{enumerate}
Therefore $\boldsymbol{g}(\cdot\,,\,\cdot)$ is a Riemmanian metric compatible
with the integrable almost-complex structure $J$ and the symplectic form $\omega$. This finishes the proof.  \end{proof}

\begin{remark}\label{metric_sobolev1/2} The metric $\boldsymbol{g}$ in the formula (\ref{kählermetric}) belongs
to the Sobolev space $H^{1/2}$.
 \end{remark}
 \begin{remark}\label{inclusionofsubgroups} We have defined loop groups
 having matrix elements of different regularity
 
 \end{remark}

\subsection{The energy functional and Birkhoff decomposition}\label{bd} Let
$f,g\in\bO$ and $\langle \,\cdot\,,\,\cdot\,\rangle$ be a left and right-invariant
metric on $U(n)$.
\begin{definition}\label{energy} The {\bf energy functional} is the map $E:\bO\to\R$ 
defined as follows:
\begin{equation}\label{energyfunctional}
E(f)=\int_{\SQ}\, \langle f'(\bu),f'(\bu)\rangle\, \bd\bu.
\end{equation}
\end{definition}
\noi As before, the inner product used is given by a Riemannian metric on $U(n)$ which is invariant under left and right translations.  Since the functions $f'$ and $g'$ belong to $H^1(\SQ,U(n))$ this function is well defined and has a gradient vector field. 

In analogy with the action of the circle group $\mathbb{T}$ on the
standard loop space (``rotating the loops'') as in \cite{AtP, Pr, Pr1, PS} 
we can ``rotate'' the adelic groups using $\SQ$. 
\begin{definition}[Rotation action]\label{rotation_flow} We recall that we are considering $\SQ$ as a multiplicative group with unit element 
$\boldsymbol1$.
The group $\SQ$ acts on $\bO$ 
as follows: if $\bv\in\SQ$, define $R_{\bv}:\bO\to\bO$ by the formula
\begin{equation}\label{rotationaction}
R_\bv(f)(\bu)=f(\bv^{-1}\bu)f(\bv^{-1})^{-1},\quad \quad f\in\bO.
\end{equation}
\end{definition} 

\noi Since based loops are preserved, 
$R_\bv(f)(\boldsymbol1)=\boldsymbol{e}$, and 
$R_{\bv_1\bv_2}=R_{\bv_1}\circ{R_{\bv_2}}$ the action is well-defined. We call this action the {\bf rotation action}. We have the self-evident:

\begin{proposition}\label{homomorphismstou1} The set of fixed points of the action of $\SQ$ on $\bO$
consists of loops $f:\SQ\to{U(n)}$  such that $f$ is a continuous (in fact real-analytic) homomorphism.

\end{proposition}

\noi We may restrict the action to the connected component of the identity
of $\SQ$. This is a one-parameter subgroup and we can use the canonical flow 
$\left\{\varphi_{_t}\right\}_{t\in\R}$ in definition (\ref{exp}).

\noi Let
$R_t:\bO\to\bO$ be defined as follows:

\begin{equation}\label{rotationflow}
R_t(f)(\bu)=f(\varphi_{_{-t}}(\boldsymbol1)\bu)(\varphi_{_{-t}}(\boldsymbol1))^{-1}.
\end{equation}

\noi The flow in equation (\ref{rotationflow}) is called the {\bf rotation flow}.

Let the Hamiltonian vector field $X_{E}$ of the energy functional $E$ be defined, as usual,
for the symplectic manifold $\bO$:
\begin{equation}\label{hamiltonianvf}
 \text{d}E(Y)=\omega(X_{E},Y),
 \end{equation}
 for any vector field $Y:\bO\to{T\bO}$ on $\bO$.
 \begin{proposition}\label{hamiltonian=rotation} The Hamiltonian flow corresponding to the energy function $E$
 is given by the rotation flow given by formula (\ref{rotationflow}). 
 In particular the rotation flow acts by symplectomorphisms.
 \end{proposition}
 
 \begin{proof} The proof is identical to the
 case of the classical loop group given in \cite{PS} Proposition (8.9.3),
 and also Lemmas 3.2, 3.3 and 3.4 in \cite{AtP}.
 Therefore, by proposition \ref{homomorphismstou1}, there exists a one-to-one correspondence between the fixed points
 of the Hamiltonian flow and the set of  continuous homomorphisms 
 $f:\SQ\to{U(n)}$.\end{proof} 
 
 \noi In fact, more generally, one has the following 
 corollary of proposition \ref{hamiltonian=rotation}:
 
 \begin{corollary} The rotation action (\ref{rotationaction}) of the group 
 $\SQ$ on $\bO$ is by
symplectomorphisms.
 \end{corollary}
 \begin{proof} It follows immediately from the fact that the connected component
 of the identity in $\SQ$ is dense, the fact that this subgroup acts by symplectomorphisms and the fact that the group of symplectomorphisms is closed in the group of diffeomorphisms with the $C^1$ topology.
 \end{proof}

 %%Proposition 2.2. The complexification? of LGL G is the loop group of the complexification of GG
%%(LG) ℂ≃L(G ℂ).

\subsection{The gradient of the energy functional}\label{energy} The 
gradient $\nabla{E}$ of the energy functional $E$, with respect to the Kähler metric $\boldsymbol{g}$ on $\bO$, is defined, as usual, as follows;

\begin{equation}\label{energygradient}
\boldsymbol{g}(\nabla E,X)=\text{d}E(X),
\end{equation}
\noi where $X:\bO\to{T\bO}$ is any vector field on $\bO$.

\begin{remark}\label{foliation-rot-grad} In a Kähler manifold the gradient and Hamiltonian vector fields
are related: the gradient is obtained by applying the complex structure to the Hamiltonian vector field. Therefore $\nabla{E}=J(X_E)$, and the vector fields $\nabla{E}$ and $X_E$ commute, since $J$ is integrable and the Nijenhuis tensor vanishes. They are the real and imaginary part of a complex vector field, hence they generate a local one-parameter complex flow, $\Psi_{T}(f)$, $T=t+\mathbf{i}s\in\mathcal{U}$, where $\mathcal{U}$ is a small neighborhood of 0 in $\C$. 

\noi Therefore, $\bO- \text{Sing}\,(\nabla{E})$ has a holomorphic 1-dimensional foliation $\mathcal{F}$ (where $\text{Sing}\,(\nabla{E})$ is the set of singularities of the gradient).  

\end{remark}
It follows the following:

\begin{corollary}\label{criticalpoints_energy} 
Since $X_E$ and $\nabla{E}$ commute, the fixed points
of the Hamiltonian action are the critical points of $E$ and vice-versa. 
Hence, the critical points of $E$ are the based loops $f:\SQ\to\bO$ which are group homomorphisms. The image of $f$ is a subgroup of $U(n)$ isomorphic to $U(1)$.
\end{corollary}

\begin{proof}
We only need to prove the last sentence. Since
$\SQ$ is compact, connected and of topological dimension one, its image must be isomorphic to the circle group $U(1)\simeq\bS^1$. \end{proof}

\noi Thus there is a one-to-one correspondence between group homomorphisms $f:\SQ\to{U(n)}$ and critical points
of $E$. 

 Let  
 \begin{equation}\label{maxtorus}
 \mathbf{T}^n=\left\{\operatorname {diag} \left(e^{i\theta _{1}},e^{i\theta _{2}},\dots ,e^{i\theta _{n}}\right):\forall j,\theta _{j}\in \mathbb {R} \right\}.
\end{equation}

\noi $\mathbf{T}^n\subset{U(n)}$ is a maximal torus of the unitary group. Let
$\mathfrak{T}^n$ be its Lie algebra (which as a real vector space is $\R^n$)
 and $\exp:\mathfrak{T}^n\to\mathbf{T}^n$ the exponential map which in this case is a homomorphism. Let $\mathbf{L}=\exp^{-1}(0)$. Then $\mathbf{L}$ is a 
 lattice in $\mathfrak{T}^n$ 
 (i.e, is a discrete, co-compact additive subgroup) so it is isomorphic to $\Z^n$.

\noi Since any homomorphism $f:\SQ\to{U(n)}$ is of the form $f=g\circ{h}$
where $h$ is a homomorphism from $\SQ$ into $U(1)$ (by corollary \ref{criticalpoints_energy} ) and $g$ is a 
homomorphism $g:U(1)\to{U(n)}$,
one has the following corollary of proposition \ref{homomorphismstou1}:

\begin{corollary}\label{coj_classes_homS1to Un} 
The set of conjugacy classes of homomorphisms from $U(1)$ to $U(n)$ is in one-to-one correspondence
with the lattice $\mathbf{L}$ up to the action of the Weyl group, which in this case is the symmetric group $S_n$ acting by permutation
of the coordinates of the lattice points. Therefore, the set of conjugacy classes of non trivial homomorphisms 
from $\SQ$ to $U(n)$ is in one-to-one correspondence with the set
$(\mathbf{L}-\boldsymbol0)/S_n\times\Q$. To the trivial homomorphism
we associate the symbol $(\boldsymbol1,0)$. 
\end{corollary}

\begin{definition}\label{cell_numbering} The group $U(n)$ acts on $\bO$ by conjugation: 
$T_{g}f(\bu)=gf(\bu)g^{-1}$. 
If $(\boldsymbol{\lambda}, q)\in(\mathbf{L}-\boldsymbol0)\times\Q$, let  
$\boldsymbol{C}_{([\boldsymbol{\lambda}], q)}\subset\bO$ denote the conjugacy class
corresponding to $(\boldsymbol{\lambda}, q)$, 
where $(\boldsymbol{\lambda}, q)$ is equivalent to $(\boldsymbol{\lambda}', q)$ if $\boldsymbol{\lambda}$ differs from $\boldsymbol{\lambda}'$ by a permutation of coordinates. 
Then, since $\boldsymbol{C}_{([\boldsymbol{\lambda}], q)}$ is an orbit of the smooth action, 
$\boldsymbol{C}_{([\boldsymbol{\lambda}], q)}$ is a finite dimensional smooth submanifold of $\bO$. In fact, $\boldsymbol{C}_{([\boldsymbol{\lambda}], q)}$ is 
diffeomorphic to the flag manifold $U(n)/\boldsymbol{T}^n$ (generically) or a quotient of it. $\boldsymbol{C}_{(\boldsymbol{1}, 0)}$ consists of only one point
corresponding to the trivial homomorphism of $\SQ$ to $U(n)$.
\end{definition}

\subsection{Action of the semigroups 
$\boldsymbol{D}^+$ and $\boldsymbol{D}^-$}
In this section we will adapt to our case the results in the papers by 
Pressley \cite{Pr, Pr1} and the book by Pressley and Segal \cite{PS}, especially 
 Theorem 8.9.9. 
 
 The punctured closed disks 
 \[\label{defD+D-}
 \boldsymbol{D}^+\overset{def}=\bar{D}^+_\Q(1)-\left\{\boldsymbol0\right\},\quad  
 \boldsymbol{D}^-\overset{def}=\bar{D}^-_\Q(1)-\left\{\boldsymbol\infty\right\},
 \]
 are abelian semigroups with identity. Let us write the elements of 
 $\boldsymbol{D}^+$ in polar coordinates as $e^{-t}\bu$ with $0\leq{t}<\infty$ and
 $\bu\in\SQ$ and the elements of $\boldsymbol{D}^-$ as 
 $e^{t}\bu$ with $0\leq{t}<\infty$.  
 Their multiplications are 
 given, respectively, by the natural formulas in polar coordinates:

\begin{equation}\label{D+}
(e^{-t_1}\bu_1)\cdot(e^{-t_2}\bu_2)=e^{-(t_1+t_2)}\bu_1\bu_2,\quad  0\leq{t_1,t_2}<\infty,\,\,\,\bu_1,\bu_2\in\SQ.
\end{equation}

\begin{equation}\label{D-}
(e^{t_1}\bu_1)\cdot(e^{t_2}\bu_2)=e^{t_1+t_2}\bu_1\bu_2, 
\quad  0\leq{t_1,t_2}<\infty,\,
\,\,\bu_1,\bu_2\in\SQ.
\end{equation}
 
\begin{proposition}[Compare, Section 7.6, 8.6 and, especially, 
Theorem 8.9.9 items (i) and (ii) of Pressley and Segal book \cite{PS} and Pressley papers \cite{Pr, Pr1}]\label{actionD+}  The negative of the gradient $-\nabla{E}$ determines a {\bf downwards gradient semi flow} (i.e. an action of the multiplicative semigroup
$\R^+=\left\{e^t\in\R : t\geq0\right\}$) $g_{e^t}:\bO\to\bO\,$ ($t\geq0$) such that:
\begin{enumerate}
\item[(\bf{i})] $g_{e^t}$ is real-analytic for all $t>0$.
\item[(\bf{ii})] $g_{e^t}$ converges to a loop which is a 
homomorphism $g_{_{\infty}}:\SQ\to{U(n)}$ as $t\to\infty$.

\end{enumerate}
The rotation action given by formula (\ref{rotationaction})
extends explicitly to an action 
$\left\{R_{e^{t}\bv}\right\}_{e^{t}\bv\in\boldsymbol{D}^+}$, 
of the semigroup $\boldsymbol{D}^+$ on $\bO:\,f\to {{R}}_{e^{t}\bv}(f)$ with
\begin{equation}\label{formula-disk-action}
{{R}}_{e^{t}\bv}(f)(\bu)=
g_{t}({R}_\bv(f)(\bu))=g_{t}(f(\bv^{-1}\bu)f(\bv^{-1})^{-1}), 
\quad t\geq0,\quad \bu,\bv \in\SQ.
\end{equation}
In the complement of the singular set of the downwards gradient, the orbits of this
action are laminations by Riemann surfaces tangent to the foliation $\mathcal{F}$
 alluded in remark (\ref{foliation-rot-grad}).
\end{proposition}
The downwards gradient semi flow $\left\{g_t\right\}_{_{t\geq0}}$ commutes with the
rotation flow $R_t$, since the vector fields $E_H$ and $\nabla{E}$ have vanishing Lie product: $[X_E,\nabla{E}]=0$. Since the one-parameter subgroup of $\SQ$, which defines the rotation flow, is dense in $\SQ$ if follows that the downwards gradient semi flow
$\left\{g_t\right\}_{_{t\geq0}}$  commutes with the rotation action:
$g_t\circ{R_v}=R_v\circ{g_t}$. 
Therefore the formula (\ref{formula-disk-action})
defines an action since 
${{R}}_{e^{t_1}\bv_1}\circ{{R}}_{e^{t_2}\bv_2}
={{R}}_{e^{t_1+t_2}\,\bv_1\bv_2}$.

\begin{remark} We will not consider the ``ascending paths" $g_t(f)$ for $t<0$.
They are defined only for $-\epsilon<t\leq0$. They are real-analytic and defined for all $t\leq0$ if and only if they are
polynomial loops (i.e the matrix coefficients are Laurent-Puiseux polynomials).  They are very interesting and important since the ``unstable manifolds"
give rise to the {\bf Bruhat decomposition}.
\end{remark}

\begin{remark}\label{descending-energy}
The energy functional is a function which is bounded below by 0. 
Hence, if we take a loop $f\in\bO$ which is not a 
critical point of $\nabla{E}$, then the energy of $g_t(f)$ is strictly 
decreasing as $t\to\infty$ i.e. it is a ``descending path". Of course a loop which is a critical point of the gradient remains fixed. It is shown that 
$g_t(f)$ converges as $t\to\infty$ to a homomorphism belonging to some manifold 
$\boldsymbol{C}_{([\boldsymbol{\lambda}], q)}$ corresponding to the conjugacy
class $([\boldsymbol{\lambda}], q)$ (see definition \ref{cell_numbering}) of $g_{_{\infty}}(f)$. The proposition implies that any $f\in\bO$ belongs to a stable manifold of some $\boldsymbol{C}_{([\boldsymbol{\lambda}], q)}$ and $\boldsymbol{C}_{([\boldsymbol{\lambda}], q)}$ is a subset of its stable manifold.
 
\end{remark}

\begin{definition}[\bf Birkhoff partition]\label{Birkhoff-partition} 
Let $S_{([\boldsymbol{\lambda}], q)}$ be the stable
manifold of $\boldsymbol{C}_{([\boldsymbol{\lambda}], q)}$ corresponding to
the conjugacy class $([\boldsymbol{\lambda}], q)$. We call
$S_{([\boldsymbol{\lambda}], q)}$ the {\bf Birkhoff stratum} with index $([\boldsymbol{\lambda}], q)$. The stable manifolds corresponding to different
indices are disjoint. It is shown in \cite{Pr1} that each $S_{([\boldsymbol{\lambda}], q)}$ is a locally closed complex submanifold, of finite codimension, of $\bO$. 
Therefore: \begin{equation}\label{stable-partition}
\bO=\underset{([\boldsymbol{\lambda}], q)}\coprod{S_{([\boldsymbol{\lambda}], q)}}
\underset{(\mathbf1,0)}\coprod{S_{(\boldsymbol1, 0)}}
\end{equation}
\end{definition}

\subsection{Filtrations of the continuous loop group} \label{filtration}
We have described adelic loops of
a Lie group $G$ having different regularity. There is a natural filtration
\begin{equation}
\La^{pol}(G)\subset\La^{rat}(G)\subset\La^{\infty}(G)
\subset\La^{\omega}(G)\subset\La^1(G)\subset\La_{W}(G)
\subset\La^2(G)\subset\La^{cont}(G)
\end{equation}\label{filtration}
\noi Here $\La^{pol}(G)$ and $\La^{rat}(G)$ are loops
such that both the loops {\it and their inverses} (so that they are groups) are given, in terms
of Fourier series, by finite Laurent-Puiseux  polynomials and series,
respectively. Some of these inclusions are homotopy equivalences but we will not discuss this here. These filtrations play an important role in the Bruhat
decomposition (see, for instance, \cite{GR}).

\section{Factorizations of adelic loop groups. Iwasawa decomposition}\label{factorizations}

 \subsection{Adelic loops as operators on Hilbert space. Adelic Grassmannian}\label{loop-operators}
 Let $\mathbb{H}(n)=L^2(\SQ,\C^n,\bd\bu)$. Then, by classical harmonic analysis, an orthonormal basis of $\mathbb{H}(n)$  is given by the functions
$\boldsymbol{\chi}^i_q:\SQ\to\C^n$ (where $1\leq{i}\leq{n},\,\, q\in\Q$) defined
by $\boldsymbol{\chi}^i_q(\bu)=\chi_q(\bu)\boldsymbol{e}_i$, where 
$\left\{\boldsymbol{e}_1,\cdots,\boldsymbol{e}_n\right\}$ is the standard
base of $\C^n$. In particular for $L^2(\SQ,M_n(\C),\bd\bu)$ we have the orthonormal
basis $\boldsymbol{\chi}^{i,j}_q(\bu)=\chi_q(\bu)\boldsymbol{e}_{i,j}$, where
$1\leq{i},j\leq{n}$ and $\boldsymbol{e}_{i,j}$ is the matrix with one in the
place $i,j$ and 0 everywhere else. Equivalently, $f\in{\mathbb{H}(n)}$ if 
$f$ has square-summable Fourier coefficients:

\begin{equation}\label{loop-in-l2}
f(\bu)=\sum_{q\in\Q}\,f_q\chi_q(\bu), \quad f_q\in\C^n,\quad \quad\sum_{q\in\Q}|f_q|^2<\infty.
\end{equation}

  Let $\text{GL}(\mathbb{H}(n))$
denote the group of invertible continuous linear operators on $\mathbb{H}(n)$. Then $\mathbb{H}(n)$
 is open in the Banach algebra $\mathcal{B}(\mathbb{H}(n))$ of bounded operators 
 of $\mathbb{H}(n)$. The group of continuous loops in $\gln$,
 $\boldsymbol\Lambda^0(\gln)$ acts as multiplication operators in 
$\mathbb{H}(n)$. 
If $\gamma:\SQ\to\gln$ and if $f\in{\mathbb{H}(n)}$ define 
$M_\gamma\in\text{GL}(\mathbb{H}(n))$ as follows:

\begin{equation}\label{loopsactionhilbertspace}
M_\gamma(f)(\bu)=\gamma(\bu)(f(\bu))
\end{equation}

\noi Then we have a continuous representation of 
$\boldsymbol\Lambda^0(\gln)$ into $\text{GL}(\mathbb{H}(n))$.

Let us decompose $\mathbb{H}(n)$ as an orthogonal direct sum of two closed subspaces
\begin{equation} \label{polarization}
\mathbb{H}(n)=\mathbb{H}^+(n)\oplus{\mathbb{H}^-(n)},
\end{equation} 
\noi where:
\[
\mathbb{H}^+(n)=\left\{f\in{\mathbb{H}(n)}\,:f(\bu)=
\sum_{q\in\Q,q\geq0} f_q\chi_q(\bu),\, f_q\in\C^n  \right\}=
\] 
\[
\left\{f\in{\mathbb{H}(n)}\,:f\,\text{is the boundary  value of a holomorphic map in}
\, D^+_\Q(1) \right\}
\]
\noi and
\[
\mathbb{H}^-(n)=(\mathbb{H}^+(n))^\perp =\left\{f\in{\mathbb{H}(n)}\,:f(\bu)
\sum_{q\in\Q,q<0} f_q\chi_q(\bu)\right\}=
\] 

\[
\left\{f\in{\mathbb{H}(n)}\,:f\,\text{is the boundary  value of a holomorphic map in}
\, D^-_\Q(1) \right\}
\]
\begin{definition}\label{canonicalpolarization} The orthogonal decomposition into two closed subspaces given in the formula \ref{polarization} is called the {\bf canonical polarization}.
\end{definition}

\begin{remark}\label{matrix-l2} If 
\[
\gamma(\bu)=\sum_{q\in\Q}\,\gamma_q\chi_q(\bu), \quad \gamma_q\in{M_n(\C)},
\]
 and we represent
the linear map $M_\gamma$ (\ref{loopsactionhilbertspace}) in terms of the Fourier basis, indexed by $\Q$, by
an infinite $\Q\times\Q$ matrix $\left(( A_{_{q_1,q_2}} \right))$  (where
$A_{_{q_1,q_2}}\in{M_n(\C)},\,\,q_1,q_2\in\Q$), then 
the associated matrix of $M_\gamma$ has elements 
$A_{_{q_1,q_2}}=\gamma_{_{q_1-q_2}}$.

\end{remark}

\medskip

The definition of the infinite  Grassmannian is the same as that given in \cite{PS}:

\begin{definition}\label{infgrass} The infinite Grassmannian $\mathbf{Gr}(\mathbb{H}(n))$
consists of all closed subspaces $W$ of $\mathbb{H}(n)$ such that:
\begin{enumerate}
\item The orthogonal projection $p_+:W\to\mathbb{H}^+(n)$ is Fredholm.
\item The orthogonal projection $p_-:W\to\mathbb{H}^-(n)$ is Hilbert–Schmidt.
\end{enumerate}
\noi $\mathbf{Gr}(\mathbb{H}(n))$ is an infinite-dimensional smooth Hilbert manifold.
\end{definition}
 
\noi Using the formula \ref{loopsactionhilbertspace} we see that
the Sobolev loop group $\La^1(\gln)$ acts linearly on $\mathbb{H}(n)$ by elements on the  Banach algebra 
$\mathcal{B}(\mathbb{H}(n))$ of bounded operators 
 of $\mathbb{H}(n)$ and as a consequence we also have a 
 continuous representation of 
$\boldsymbol\Lambda^1(\gln)$ into $\text{GL}(\mathbb{H}(n))$.

\begin{definition}
Define $\text{GL}_{\text{R}}(n)\subset\text{GL}(\mathbb{H}(n))$ 
as a matrix of bounded operators given by blocks with respect to the polarization (\ref{polarization}),
as follows:

\begin{equation}\label{restrictedglH}
\text{GL}_{\text{R}}(n)=
\left\{
\begin{pmatrix}
a & b \\
c & d 
\end{pmatrix} 
\,:\, a,\, d \,\,\text{are Fredholm and} \,\,
c ,\, d\,\,\text{are Hilbert-Schmidt}
\right\}.
\end{equation}
\noi ${GL}_R(n)$ is called 
the {\bf restricted linear group}. The group 
${\mathbb{U}}_{R}(n)={\mathbb{U}(n)}\cap\text{GL}_{\text{R}}(n)$, where
$\mathbf{\mathbb{U}}(n)$ is the unitary group of $\mathbb{H}(n)$ is called 
the {\bf restricted unitary group}.
\end{definition}

\begin{proposition} 
If $W\in\mathbf{Gr}(\mathbb{H}(n))$ and $F\in\text{GL}_{\text{R}}(n)$ then 
$F(W)\in\mathbf{Gr}(\mathbb{H}(n))$. Therefore, $\text{GL}_{\text{R}}(n)$
and $\mathbb{U}_R(n)$ act smoothly on the restricted Grassmanian. 
Furthermore, the action of both groups is transitive.
\end{proposition}
\noi The proof of this proposition is in 
\cite{PS} Chapter 7 but is also given in full detail in \cite{Lau}.
The proof of the transitivity of the action of the restricted unitary group 
$\mathbb{U}_R(n)$ on the restricted Grassmannian is essentially an application 
of the Gram–Schmidt process.

\begin{proposition}\label{actionongrassmanian}
The group $\La^1(\gln)$ acts by formula (\ref{loopsactionhilbertspace}) on  
$\cH(n)$. If $\gamma$ belongs to $\La^1(\gln)$ and 
$M_\gamma\in\text{GL}(\calH_1(n))$ is defined by formula (\ref{loopsactionhilbertspace}), then its block decomposition (\ref{restrictedglH}), with respect to the
canonical polarization, satisfies the conditions to belong to the restricted linear group $\text{GL}_{\text{R}}(n)$. Therefore $\La^1(\gln)$ acts on 
$\mathbf{Gr}(\mathbb{H}(n))$. In particular 
the loop group $\bO$ acts on $\mathbf{Gr}(\mathbb{H}(n))$. 
\end{proposition}

\begin{proof} The proof is similar to that of Propositions 2.3 and 2.7 in \cite{SW}. If $\gamma\in\La^1(\gln)$ has Fourier expansion
\begin{equation}\label{fs}
\gamma(\bu)=\sum_{q\in\Q}\,\gamma_q\chi_q(\bu), \quad \gamma_q\in{M_n(\C)}
\end{equation}
\noi then, necessarily:
\[
\sum_{q\in\Q}|\gamma_q|<\infty,
\]
\noi because functions which belong to Sobolev space $H^1(\SQ,M_n(\C))$ also
belong to the Wiener algebra. In terms of the countable orthonormal basis of
 of $\mathbb{H}(n)$, the linear map $M_\gamma$ can be expressed
 as a $\Q\times\Q$ infinite matrix given in remark \ref{matrix-l2} and it can be verified, using the fact that the series (\ref{fs}) is absolutely convergent, that the block
 matrices $a,b,c,d$ satisfy the required conditions. \end{proof}

\noi As in the case of the Wiener algebra, using the fact that 
$H^1(\SQ,\C^n)$ is a Banach algebra we define  the Banach subalgebras
 $\La_+^1(\gln)$ and  $\La_-^1(\gln)$ as follows:
 
 \begin{definition}\label{splittinglambda+-}
Let $\La_+^1(\gln)$ denote the set of functions in
 $\La^1(\gln)$ such that
 $f$ extends to a continuous map $f:\bar{D}^+_\Q(1)\to\gln$ which is holomorphic
 in ${D}^+_\Q(1)$. Analogously,
 let $\La_-^1(\gln)$ denote the set of functions in
 $\La^1(\gln)$ such that
 $f$ extends to a continuous map $f:\bar{D}^-_\Q(1)\to\gln$ which is holomorphic
 in ${D}^-_\Q(1)$.

 \end{definition}
 \subsection{Sobolev Grassmannian model for the adelic loop group $\bO$}\label{sobgrass}
 Let $\mathcal{H}(n)= H^1(\SQ,\C^n,\bd\bu)$ be the Hilbert space of functions
 in $H^1(\SQ,\C^n)$ with inner product given by:
 \begin{equation}
 \langle\, f \,,\, \,g\rangle_{_1}=
 \int_{\SQ}\,\langle\, f(\bu) \,,\, \,g(\bu)\rangle_{_{\C^n}}\bd\bu+
 \int_{\SQ}\,\langle\, f'(\bu) \,,\, \,g'(\bu)\rangle_{_{\C^n}}\bd\bu.
 \end{equation}\label{innproduct}
 
\noi In  formula $\langle\,\cdot\,,\,\cdot\rangle_{_{\C^n}}$ is the standard hermitian product of $\C^n$ and $f'$, $g'$ are weak derivatives.
 In terms of Fourier series: 
 \begin{equation}
 \left<\, \sum_{q\in\Q}a_q\chi_q \,,\, \,\sum_{q\in\Q}b_q\chi_q\right>_{_1}=
 \sum_{q\in\Q}\,\langle\,a_q\,,\,b_q\rangle_{_{\C^n}}(1+|q|^2)
 \end{equation}
 Let $\calH_1(n)=\calH_1^+(n)\oplus\calH_1^-(n)$ be the polarization of 
 $\calH_1(n)$ where
  $\calH_1^+(n)$ and
 $\calH_1^-(n)$ are defined exactly as in the canonical polarization.
 Define the restricted Sobolev Grassmannian $\mathbf{Gr}_1(n)$ as the set of closed subspaces
 $W$ of $\calH_1(n)$ satisfying the same conditions of 
 definition \ref{infgrass}. Analogously we define the restricted Sobolev group
 of automorphisms as follows:

 \begin{definition}\label{sobolev-gln} Let $\text{GL}_1\calH(n)$ be the group of linear automorphisms of $\calH_1(n)$ which can be decomposed in blocks, with respect to the polarization, as the matrix
 in (\ref{restrictedglH}) where $a,b,c, d$  satisfy the same conditions:
 $a,d$ Fredholm and $b,c$ Hilbert-Schmidt. $\text{GL}_1\calH_1(n)$ is called
 the {\bf restricted Sobolev linear group}.
 \end{definition} 
 
 As pointed in \cite{PS} section 7.2, there are several dense submanifolds
 of $\mathbf{Gr}(\mathbb{H}(n))$ which are Grassmannians
 corresponding to the degree of regularity
 of the functions in the subspaces. An important role there is the case of
 the smooth Grassmannian. For us the relevant Grassmannian 
 is $\mathbf{Gr}_1(n)$.
 
 \medskip
 \noi The following proposition is an immediate consequence of the fact that
 $H^1(\SQ,\C,\bd\bu)$ is dense in $L^2(\SQ,\C,\bd\bu)$.
 
\begin{proposition}\label{inclusionhilbertinsobolev} The natural inclusion
$\calH_1(n)\subset\cH(n)$ is compatible with the canonical polarization:
$\calH_1^+(n)\subset\cH^+(n)$ and $\calH_1^-(n)\subset\cH^-(n)$. Furthermore,
 $\calH_1(n)$ is dense in $\cH(n)$. Therefore, there is a natural
 inclusion of Grassmannians $\mathbf{Gr}_1(n)\subset\mathbf{Gr}(\mathbb{H}(n))$
 and $\mathbf{Gr}_1(n)$ is dense in $\mathbf{Gr}(\mathbb{H}(n))$. In addition,
 since $H^1(\SQ,\C,\bd\bu)$ is closed under multiplication by proposition
 \ref{h1Banach-algebra} it follows that $\La^1(\gln)$ 
 acting on $\mathbf{Gr}(\mathbb{H}(n))$ leaves invariant $\mathbf{Gr}_1(n)$.
 \end{proposition}
 
 \noi The next proposition can be proven exactly as the corresponding
 proposition \ref{actionongrassmanian}:

 \begin{proposition}
 The group $\La^1(\gln)$ acts by formula (\ref{loopsactionhilbertspace}) on  
$\calH_1(n)$. If $\gamma$ belongs to $\La^1(\gln)$ and 
$M_\gamma\in\text{GL}(\calH_1)$
has the block decomposition (\ref{restrictedglH}), with respect to the
canonical polarization of $\calH_1(n)$, then $a,b,c,d$ satisfy the conditions in \ref{restrictedglH} to belong to $\text{GL}_1\calH_1(n)$. Therefore $\La^1(\gln)$ acts on the restricted Grassmannian
$\mathbf{Gr}_1(n)$. In particular 
the loop group $\bO$ acts on $\mathbf{Gr}_1(n)$. 

\noi The group 
${\mathbb{U}}_1(n)={\mathbb{U}_1(\calH_1(n))}\cap\text{GL}_1\calH_1(n)$, where
${\mathbb{U}_1(\calH_1(n))}$ is the unitary group of $\calH_1(n)$, is called 
the {\bf restricted Sobolev unitary group}.
\end{proposition}

 For every rational $q\geq0$ the orbit of $\,\calH_1^+(n)$ under 
 $\La^1(\gln)$ contains 
 $\chi_q\calH_1^+(n)$, where this space consists of the functions
 in $\calH_1^+(n)$ multiplied by the scalar function given
 by the character  $\chi_q\,\,(q\geq0)$ of $\SQ$. 
 This is because multiplication by a character commutes with 
 the action of $\La^1(\gln)$. This motivates the following definition:

\begin{definition}

\noi $\boldsymbol{Gr}_\Q^{(n)}$ 
denotes the subset of closed subspaces $W\in\mathbf{Gr}_1(n)$
such that:

\noi $\chi_{_{1/n}}\,W\subset{W},\,\,\forall\,n\,\in\N$. Where 
$\chi_{_{1/n}}W$ consists of the product of the functions in 
$H^1(\SQ,\C^n)$ belonging to $W$ by the scalar function determined 
by the character $\chi_{_{1/n}}$ of $\SQ$. 
Obviously $\calH_1^+(n)\in\boldsymbol{Gr}_\Q^{(n)}\bf{.}\,\,$
{\bf $\boldsymbol{Gr}_\Q^{(n)}$ is called the Sobolev adelic Grassmannian model of $\bO$.}
\end{definition}

\begin{remark}In the case of the standard loop group the orbit of 
 $\mathbb{H}^+(n)$ is characterized as the set of subspaces $W$ 
 such that $zW\subset{W}$, i.e., one uses only one character of $\Z$ (see section 8.3 in \cite{PS}). The definition is also motivated by the fact that $\Q$ is the direct limit of the infinite cyclic groups generated by $1/n,\,n\in\N$:
$\Q\cong \underset{\underset{n\in\N}{\longrightarrow} }\lim\,\frac{1}{n}\Z$.
 \end{remark}

\begin{remark}\label{W+} If $W\in\boldsymbol{Gr}_\Q^{(n)}$ let 
 $W^+=\overline{\text{Span}\left\{\chi_{_{1/n}}W: n\in\N\right\}}$ i.e., the
 closed subspace of $W$ which is the closure of the space generated by the union of the subspaces $\chi_{_{1/n}}W$. Exactly the same argument used in proof 8.3.2
of \cite{PS} page 126, implies that $W^+$ is a closed subspace of $W$ of codimension $n$.
\end{remark}
 
 \begin{theorem}\label{grass=loop} The Sobolev loop group $\La^1(U(n))$ of the unitary group
 $U(n)$, acts transitively on $\boldsymbol{Gr}_\Q^{(n)}$. The isotropy subgroup
of $\calH_1^+(n)$ consists of constant loops
$\gamma:\SQ\to{U(n)}$.
 \end{theorem}
 \begin{proof} In this proof we regard $\boldsymbol{Gr}_\Q^{(n)}$ embedded
 in the $L^2$ Grassmannian $\mathbf{Gr}(\mathbb{H}(n))$.
 The proof is similar to that of \cite{PS} so we only indicate how to adapt it to our case. Let $W\in\boldsymbol{Gr}_\Q^{(n)}$ and let $W^P$ be the orthogonal complement (with respect to the $L^2$ inner product of $\cH(n)$) of $W^+$ in $W$.
 
\noi Let $w_1,\cdots,w_n$ be an orthonormal basis  (again with respect to
 the $L^2$ product) of $W^P$.

 \noi If in terms of Fourier series $w_k(\bu)=\sum_{q\in\Q}a_{_{k,q}}\chi_q(\bu),\,\, a_{_{k,q}}\in\C^n$ for $k=1,\cdots,n$, then:
 
 \begin{equation}\label{unitarygamma}
 \begin{split}
 \langle\, w_i(\bu) \,,\, \,w_j(\bu)\rangle_{_{\C^n}} & = \sum_{p,q\in\Q}\,\langle\,a_{_{i,q}}\chi_q(\bu)\,,\,a_{_{j,p}}\chi_p(\bu)\rangle_{_{\C^n}} \\
 & = \sum_{p,q\in\Q}\,\langle\,a_{_{i,q}}\,,
 \,a_{_{j,p}}\rangle_{_{\C^n}}\chi_{_{q-p}}(\bu)\\
 & = \sum_{q\in\Q}\,\langle\,\chi_qw_i\,,
 \,w_j\rangle_{_{L^2}}\chi_{_{q}}(\bu)=\delta_{ij}\,, \quad\quad 1\leq i, j, \leq n.
  \end{split}
 \end{equation}
 \noi Therefore, if $\mathbf{w}(\bu)=\left(w_i(\bu)  \right)$ is the matrix with
 columns $w_i$ we have by equation (\ref{unitarygamma}), that 
 $\mathbf{w}:\SQ\to{U(n)}$ and hence $\mathbf{w}\in\La^1(U(n))$. The space
 $(\calH_1^+(n))^P$ consists only of constant functions from $\SQ$ to $\C^n$. It
 follows that $M_\mathbf{w}((\calH_1^+(n))^P)=W^P$ since $M_\mathbf{w}(\mathbf{e_i})=w_i$ 
 (here $\mathbf{e_i}$ is a column vector in the canonical basis of $\C^n$)
 and this implies that $M_\mathbf{w}(\calH_1^+(n))=W$, because
 $M_\mathbf{w}((\calH_1^+(n) +\chi_{1/n}\calH_1^+(n))^P)
 =M_\mathbf{w}((\calH_1^+(n))^P)+M_\mathbf{w}((\chi_{1/n}\calH_1^+(n))^P)
 \subset{W},\,\,\forall\,n\in\N$ and
 $\underset{n\in\N}\cap{\chi_{1/n}\calH_1^+(n)}=\mathbf0$. Hence
 the action of $\La^1(U(n))$ on $\boldsymbol{Gr}_\Q^{(n)}$ is transitive.
 This way
 we have associated to each subspace $W\in\boldsymbol{Gr}_\Q^{(n)}$ a loop 
 $\mathbf{w}\in\La^1(U(n))$. 
 
 \medskip
 \noi To finish the proof we only have to compute the isotropy
 subgroups of the action. It is clear that the stabilizer of $\calH_1^+(n)$ of the action
 of $\La^1(\gln)$ is
 $\La_+^1(\gln)$ (definition (\ref{splittinglambda+-})). Hence the stabilizer
 is a loop $\gamma$ in $\La_+^1(U(n))\cap\La_+^1(\gln)$.
 However we claim that every loop $\gamma\in\La_+^1(U(n))\cap\La_+^1(\gln)$
 is constant. Indeed, $\gamma$ can be extended to a map
 (denoted with the same symbol): $\gamma:\bar{D}^+_\Q(1)\to{M_n(\C)}$.
 Consider the map $T$, given by the square of the Frobenius norm of the matrix $\gamma(\bu)$: 
  \[
 T(\bz)=|\gamma(\bz)|_{_{Fr}}=
 |\gamma(t\bu)|_{_{Fr}}=\text{Trace}
 (\gamma^*(t\bu)\gamma(t\bu)),\quad \bz=t\bu;\,\, \bu\in\SQ,\,\,0\leq{t}\leq1,
 \] 
 
 \noi where $\gamma^*(\bu)$ is the  conjugate transpose of $\gamma(\bu)$. 
 Then the map $T$ is constant and equal to $n$ in $\SQ$,  since  
 $\gamma(\bu)\in{U(n)}$. On the other hand $\gamma$ is holomorphic
 on the open disk ${D}^+_\Q(1)$; then by the maximum principle 
 for matrix-valued functions using the Frobenius norm 
 (see Theorem 1, {\it Maximum Frobenius norm principle} in \cite{Co}), which is a consequence of our proposition \ref{Liouville}, 
 $\gamma(\bu)=U_0,\,\,\forall\bu\in\SQ$, for some $U_0\in{U(n)}$ i.e., $\gamma$ is a constat loop. \end{proof} 
 \begin{corollary} The adelic Grassmannian
 $\boldsymbol{Gr}_\Q^{(n)}$ is diffeomorphic to 
 $\La^1(U(n))/U(n)$. Since $\La^1(U(n))/U(n)$ is canonically 
 identified with the based loop group $\bO$.  This justifies the name
 {\bf Sobolev adelic Grassmannian model of $\bO$.} 

  \end{corollary}
  \begin{corollary}\label{grass-loops-lambdas}
  $\boldsymbol{Gr}_\Q^{(n)}$ is
  an infinite-dimensional Lie group which is a Kähler manifold. The action on itself by left translations is by biholomorphisms.
  On the other hand $\La^1(\gln)$ also acts transitively on $\mathbf{Gr}_1(n)$. The isotropy subgroup of $\calH_1^+(n)$ is $\La_+^1(\gln)$. 
  
  \noi Therefore we have the following equality of Kähler homogeneous manifolds:
  
\begin{equation}\label{grass-loops-lambdas-eq}
\bO= \La^1(\gln)/\La^1_+(\gln)=\La^1(U(n))/U(n)=\boldsymbol{Gr}_\Q^{(n)}.
\end{equation}

 \end{corollary}
 
 \medskip 
The proof of theorem \ref{grass=loop} actually proves that
any loop $\gamma\in\La^1(\gln)$ factorizes as 
$\gamma=\mathbf{w}\cdot{\gamma_+}$ with $\mathbf{w}\in\La^1(U(n))$ and
$\gamma_+\in\La^1(\gln)$. If $\gamma=\mathbf{w}'\cdot{\gamma_+^{\prime}}$
with $\mathbf{w}'\in\La^1(U(n))$ and
$\gamma_+^{\prime}\in\La^1(\gln)$ is other factorization then
$\mathbf{w}'\cdot\mathbf{w}^{-1}$ belongs to $\La^1(U(n))\cap\La^1(\gln)$ and,
as we showed, the Frobenius maximum modulus principle implies that
$\mathbf{w}'\cdot\mathbf{w}^{-1}$ is a constant loop in $\La^1(U(n))$.

  This leads to the following theorem:

 \begin{theorem}[\bf Iwasawa decomposition]\label{iwasawadecomposition} The multiplication map:
 \begin{equation}
 F:\bO\times\La^1_+(\gln)\to\La^1(\gln)
 \end{equation}
 given by  $(f,g)\mapsto{f\cdot{g}}$ 
 where $({f}\cdot{g})(\bu)=f(\bu)g(\bu)$, $\bu\in\SQ$,
 is a diffeomorphism onto. In particular, any loop $\gamma\in\La^1(\gln)$
 has a unique decomposition $\gamma=\mathbf{w}\cdot\gamma_+$ where $\mathbf{w}\in\bO$ and
 $\gamma_+\in\La^1_+(\gln)$.
 \end{theorem}
\begin{proof} The proof follows easily from theorem \ref{grass=loop} and its corollaries because one can verify that the
projection between Banach manifolds $\Pi:\La^1(\gln)\to\bO$ is a locally trivial differentiable
fibration (in fact a principal fibration) and by the existence and uniqueness of the decomposition $\gamma=\mathbf{w}\cdot\gamma_+$. This proves that $\Pi$ is a trivial principal fibration.
Proofs of this theorem for the standard loop group of $U(n)$ can be found, for instance, in \cite{BD, Ke, PS} \end{proof}

\section{Birkhoff Factorization}\label{BF} 

The purpose of this section is to state and prove the following theorem:

 \begin{theorem}[{\bf Birkhoff factorization}] \label{BFT}\hfil
 \begin{enumerate}
\item Any loop
$\gamma:\SQ\to\gln$, $\gamma\in\La^1(\gln)$, can be factorized
as follows:
\begin{equation}\label{-d+}
\gamma=\gamma_-\cdot\Delta_{_{\mathbf{q}}}\cdot\gamma_+,
\end{equation}
\noi $\gamma_\pm\in\La_\pm^1(\gln)$, 
$\mathbf{q}=(q_1,\cdots,q_n)\,\, (q_i\in\Q\,\,\,\forall{i})$,
$\Delta_{_{\mathbf{q}}}:\SQ\to{U(n)}$ is the diagonal homomorphism:
 \begin{equation}
 \Delta_{_{\mathbf{q}}}(\bu)=
 \begin{pmatrix}
    \chi_{q_1}(\bu) & & \\
    & \ddots & \\
    & & \chi_{q_n}(\bu)
  \end{pmatrix}
 \end{equation}
 
 \item In any other factorization of the type given by equation (\ref{-d+})
 the diagonal middle factor differs from $\Delta_{_{\mathbf{q}}}$ only
 in the order of the diagonal elements.

 \end{enumerate}
 \end{theorem}
\begin{proof}
We recall that $\bO$ admits
the partition (\ref{stable-partition}). The set
$S_{([\boldsymbol{\lambda}], q)}$ is the stable manifold of the conjugacy
class $\boldsymbol{C}_{([\boldsymbol{\lambda}], q)}$
of the homomorphism
$\lambda:\SQ\to\bO$, with respect to the descending gradient flow 
of the energy functional. In \cite{Pr1} it is shown that 
$S_{([\boldsymbol{\lambda}], q)}$ is a locally closed complex submanifold
of $\bO$. 
Theorem \ref{iwasawadecomposition} implies that
the projection map $\Pi:\La^1(\gln))\to\bO$ is a principal fibration
and we obtain, using the partition \ref{stable-partition}, the partition: 

\begin{equation}
\La^1(\gln))=\underset{([\boldsymbol{\lambda}], q)}
\coprod{\Sigma_{([\boldsymbol{\lambda}], q)}}\coprod{S_{(\boldsymbol1, 0)}}
\end{equation}
\begin{equation}
\Sigma_{([\boldsymbol{\lambda}], q)}=\Pi^{-1}\left\{S_{(\boldsymbol{\lambda}], q)}
\right\} \,\, \text{and} \,\, 
\Sigma_{(\boldsymbol1, 0)}=\Pi^{-1}\left\{S_{(\boldsymbol1, 0)}\right\}
\end{equation}

\noi Each stratum $\Sigma_{([\boldsymbol{\lambda}], q)}$ contains
the finite dimensional submanifold  
$\boldsymbol{C}_{([\boldsymbol{\lambda}], q)}
\subset{U(n)}\subset\La^1(\gln),$ 
of definition \ref{cell_numbering}, consisting of the conjugacy class
of the homomorphism $\lambda:\SQ\to{U(n)}$.

On the other hand, by corollary \ref{grass-loops-lambdas} the loop
group is the homogeneous space of right-cosets
$\bO= \La^1(\gln)/\La^1_+(\gln)$, therefore the subgroup 
$\La^1_-(\gln)$, 
consisting of the loops 
$\gamma:\SQ\to\gln$ which can be extended holomorphically to $D^-_\Q$,
 acts on the left on the space of right-cosets. 
 
 Specifically: 
 by theorem \ref{iwasawadecomposition} evey right-coset can be uniquely
 represented as the right-coset $\mathbf{w}[\La^1_+(\gln)]$ with
 $\mathbf{w}\in\bO$.  Let $[\mathbf{w}]$
denote this coset. If $\gamma_{_{-}}\in\La^1_-(\gln)$ and 
$(\gamma_{_{-}}\cdot\mathbf{w})[\La^1_+(\gln)]$ is the right-coset obtained by 
left translating,
by $\gamma_{_{-}}$ the coset $[\mathbf{w}]$, then this right-coset is of the form  
$\mathbf{w}'[\La^1_+(\gln)]$ for a unique $\mathbf{w}'\in\bO$. The left action of 
$\La^1_-(\gln)$ is given by
$L_{\gamma_{_{-}}}([\mathbf{w}])=[\mathbf{w}'],\,\,\,\gamma_{_{-}}\in\La^1_-(\gln)$.

Consider the orbit $\mathcal{O}(\bw)$ of $\bw$ under the action of the semigroup
$\boldsymbol{D}^+$ given by the formula (\ref{formula-disk-action}) 
in proposition (\ref{actionD+}):
$$
\mathcal{O}(\bw)=\left\{{{R}}_{e^{t}\bv}(\bw):t\geq0,\, \bv\in\SQ\right\}
$$
The parametrization of the orbit gives a map 
$F_{\bw}:\boldsymbol{D}^+\to\bO,\quad 
t\bv\mapsto{R}_{e^{t}\bv}(\bw)$.

By Formula 7.6.2 of Proposition 7.6.1 in \cite{PS} this orbit can be extended by adding the point 
$\underset{t\to\infty}\lim{{R}}_{e^{t}\mathbf1}(\bw)=
\underset{t\to\infty}\lim{g_t(\bw)}$ which is independent of $\bw$.
This implies that this parametrization $F_{\bw}$
extends to a continuous map $\bar{F}_{\bw}:\bar{D}_\SQ^+\to\bO$ by
defining $\bar{F}(\boldsymbol0)=\underset{t\to\infty}\lim{{R}}_{e^{t}\mathbf1}(\bw)=\underset{t\to\infty}\lim{g_t(\bw)}\overset{def}=\lambda$. By proposition
\ref{actionD+} item {\bf(ii)}, $\lambda$ is a group homomorphism  
$\lambda:\SQ\to{U(n)}\subset\gln$, therefore it belongs to a unique stratum
$S_{([\boldsymbol{\lambda}], q)}$. 

\noi If 
$\,\overline{\mathcal{O}(\bw)}=F_{\bw}(\bar{D}^+_\SQ)$ then, since every point
of this set is in the stable manifold of 
$\boldsymbol{C}_{([\boldsymbol{\lambda}], q)}$,  it follows that 
$\,\overline{\mathcal{O}(\bw)}\subset{S_{([\boldsymbol{\lambda}], q)}}$.
Therefore, the union of the one-parameter family of right-cosets
$[g_t(\bu)],\,t\geq0\,,$ is contained in $\Sigma_{([\boldsymbol{\lambda}], q)}$. 
In particular the right-coset
$[\lambda]=\lambda[\La^1_+(\gln)]$ is contained in 
$\Sigma_{([\boldsymbol{\lambda}], q)}$.  

We now will prove that that the stratum
$\Sigma_{([\boldsymbol{\lambda}], q)}$ belonging to the conjugate class
of $\lambda$ is equal to the orbit of the set
of right-coset $[\lambda]$ under the left action of $\La^1_-(\gln)$.
Let $D_\Q^+(1)$ be the semigroup with product
in polar coordinates: ${t_1\bv_1\cdot{t_2}\bv_2}=t_1t_2\bv_1\bv_2$.
Define the action of this semigroup as the 
``{\it scaling semi flow}'':
\[
\varphi_{_{t\bv}}:\La^1_-(\gln)\to\La^1_-(\gln),\quad {t}\bv\in{D_\Q^+(1)},
\]

\noi which is  given in terms of Fourier series as follows:

\begin{equation}
\varphi_{_{t\bv}}(\gamma_{_{-}})(\bu)=\sum_{q\leq0}t^q\,a_q\chi_q(\bv\bu)
=\sum_{q\leq0}t^q\,a_q\chi_q(\bv)\chi_q(\bu),\quad 0\leq{t}\leq1,\,\, \bu, \bv\in\SQ,
\end{equation}
\noi if $\gamma_{_{-}}\in\La^1_-(\gln)$  is given by the Fourier series 
$\gamma_{_{-}}(\bu)=\underset{q\leq0}\sum \,a_q\chi_q(\bu)$. Clearly
$\varphi_{_{t_1\bv_1\cdot{t_2}\bv_2}}=
\varphi_{_{t_1\bv_1}}\circ\varphi_{_{t_2\bv_2}}$ and 
$\varphi_{_{\boldsymbol0}}(\gamma_{_{-}})$ is the constant loop 
$\gamma_{_{-}}(\boldsymbol\infty)=a_{_0}\in\gln$. Hence 
the one-parameter family of right-cosets 
$(\varphi_{_{t\bv}}(\gamma_{_{-}})\cdot\lambda)[\La_+^1(\gln)]$ belongs
to the orbit of the left action of $\La^1_-(\gln)$ $[\lambda]$. In particular
when $t=0$ we have $a_{_0}\lambda[\La_+^1(\gln)]=
a_{_0}\lambda{a^{-1}_{0}}[\La_+^1(\gln)]$, and $a_{_0}\lambda{a^{-1}_{0}}$ is a diagonal matrix which differs from $\lambda$ only by a permutation of the
diagonal elements.
Furthermore $\Sigma_{([\boldsymbol{\lambda}], q)}$ is the stable manifold
of $\boldsymbol{C}_{([\boldsymbol{\lambda}], q)}$
and thus the right-coset 
$[\lambda]$ belongs to the orbit of $\La^1_-(\gln)$. We summarize what we have shown:
\begin{equation}\label{BirkhoffCells}
 \La^1(\gln)= 
 \end{equation}\label{}
\hfil $\underset{([\boldsymbol{\lambda}], q),q\neq0}\coprod\,
\La_-^1(\gln)\cdot\boldsymbol{C}_{([\boldsymbol{\lambda}], q)}
\cdot\La_+^1(\gln)\coprod\La_-^1(\gln)\cdot \La_+^1(\gln)$,

\noi with the obvious meaning of the triple products. 
This obviously completes the proof of theorem \ref{BF}. \end{proof}

\begin{remark} Because of the equality of the homogeneous spaces in
equation (\ref{grass-loops-lambdas-eq}) we could have done in the previous theorem
using the adelic Grassmannian $\boldsymbol{Gr}_\Q^{(n)}$.
\end{remark}

\begin{remark} The Birkhoff factorization for the
classical loop groups is determined by points in a an integer lattice
modulo permutation of the coordinates but in the case of adelic groups
there is a ``weight'' $q\in\Q$ except for the trivial loop.
That is why there is the distinguished component $\La_-^1(\gln)\cdot \La_+^1(\gln)$.
\end{remark}
\begin{definition} The elements 
$B([\boldsymbol{\lambda}], q)\overset{def}=\La_-^1(\gln)\cdot
\boldsymbol{C}_{([\boldsymbol{\lambda}], q)}\cdot\La_+^1(\gln)$, $q\neq0$, 
and $B(\mathbf1, 0)\overset{def}=\La_-^1(\gln)\cdot \La_+^1(\gln)$
of the partition (\ref{BirkhoffCells}) are called {\bf Birkhoff cells}.
They are indexed by the conjugacy classes of homomorphisms
$f:\SQ\to{U(n)}$.
\end{definition}

\begin{corollary}[\bf The Big-cell] The use of the scaling flow shows that 
the the Birkhoff cell $B([\boldsymbol{\lambda}], q)$ retracts strongly to 
the manifold $\boldsymbol{C}_{([\boldsymbol{\lambda}], q)}$ if $q\neq0$. If
$q=0$ the set $\boldsymbol{C}_{(\mathbf1, 0)}$ consists of a single point
which is the constant loop equal to the identity matrix $I$ in $\gln$ and the
Birkhoff cell $B(\mathbf1, 0)$ retracts strongly to this point so that 
 $B(\mathbf1, 0)$ is contractible. In fact $B(\mathbf1, 0)$ is open and
 dense in $\La^1(\gln)$. If $\La_{-,*}^1(\gln)$ is the subgroup
 of $\La_-^1(\gln)$ consisting of loops $\gamma_{_{-}}$ such that 
 $\gamma_{_{-}}(\mathbf1)=\text{I}$. Then the map:
 \[
 F:\La_{-,*}^1(\gln)\times\La_+^1(\gln)\to\La^1(\gln)
 \]
 \noi is a diffeomorphism onto the open and dense Birkhoff cell $B(\mathbf1, 0)$.
 For this reason $B(\mathbf1, 0)$ is called the {\bf big cell}. This implies  that any loop $\gamma$ in the big cell can be written uniquely
 as a product $\gamma=\gamma_{_{-}}\cdot\gamma_{_{+}}$. We remark that
 $\La_{-,*}^1(\gln)\cap\La_+^1(\gln)=\left\{I\right\}$ because the loops
 in the intersection can be extended holomorphically to $\PQ$, so by Liouville's theorem \ref{Liouville} they are constant loops. The map $F$ is differentiable because $\La^1(\gln)$ is a Banach Lie group. One can verify directly that
 the differential at every point determines an isomorphism of the Banach algebra
 $\mathcal{L}(\gln)$. The big cell is clearly open and dense.
 
 \end{corollary}
 
 \section{Birkhoff-Grothendieck theorem for $\PQ$}\label{BG}
 
\begin{theorem}[\bf Birkhoff-Grothendieck theorem for $\PQ$] 
\cite{Bir, Bir1, Grot}.
 Given a holomorphic vector bundle $\pi:E\to\PQ$ of rank $n$ over $\PQ$, there exists $\boldsymbol{q}=(q_q,\cdots,q_n),\,\,q_i\in\Q$, such that $E$ is holomorphically isomorphic to the direct sum of line bundles:
  
 \begin{equation}\label{BGdecomposition}
 E=\I_{\C P^1_\Q}(q_1)\oplus\cdots\oplus\I_{\C P^1_\Q}(q_n),
 \end{equation} 
 
\noi where $\I_{\C P^1_\Q}(q_i)$ is the pullback of the line 
bundle $\mathcal{O}_{\C P^1}(1)$ of Chern class 1
 over $\C P^1$ by a map $f:\PQ\to{\C P^1}$ (see section (\ref{vb}))
 in the homotopy class in $[\PQ,\C P^1]$  corresponding to the
 rational number $q_i$. Any other decomposition of this form differs only
 by a permutation of the summands.
 \noi
 \end{theorem}
 
\begin{proof} The proof that the Birkhoff factorization theorem
implies the theorem is standard. By definition, the bundle $E$ 
is given by a single holomorphic cocycle:
\[
f_{_{(r,R)}}(z): A(r,R)\to{\text{GL}}(n,\C).
\]
We can assume without loss of generality that $r<1$ and $R>1$ (because
for any $s>0$ the map $t\bu\mapsto{st\bu}$ is holomorphic). Then the annulus
$A(r,R)$ contains $\SQ$ and the restriction of $f_{_{(r,R)}}$ to $\SQ$ is a loop
$\gamma:\SQ\to\gln$. By theorem (\ref{BFT}), 
$\gamma=\gamma_{_{-}}\cdot\Delta_{\mathbf{q}}\cdot\gamma_{_{+}}$ for some
$\mathbf{q}=(q_1,\cdots,q_n)$. The fact that
$\gamma$ extends as the function $f_{_{(r,R)}}$ implies that both loops
$\gamma_{_{-}}$ and $\gamma_{_{+}}$ extend also as holomorphic maps
on the annulus. But then the cocycle $f_{_{(r,R)}}$ is holomorphically equivalent
to the cocycle 
$\hat\Delta_{_{\mathbf{q}}}:A(r,R)\to\gln,\,\,$ 
$\hat\Delta_{_{\mathbf{q}}}(t\bu)=
\text{diagonal}\,(t^{q_1}\chi_{q_1}(\bu),\cdots,t^{q_n}\chi_{q_n}(\bu))\,$ 
($r<t<R$),
and this holomorphic cocycle defines the direct sum
in formula (\ref{BGdecomposition}).
That $\gamma_{_{-}}$ and $\gamma_{_{+}}$ extend as holomorphic maps
on the annulus is because 
$\gamma_{_{-}}\cdot\hat\Delta_{_{\mathbf{q}}}$ and 
$f_{_{(r,R)}}(z)\cdot\gamma_{_{-}}^{-1}$ are holomorphic in the annuli
$A(r,1)$ and $A(1,R)$, respectively,  and they coincide on $\SQ$ and since the Laurent coefficients of a holomorphic
function on the annulus $A(r,R)$ are given by a Cauchy-type integral
\begin{equation}
a_{q}=\int_{\SQ}\chi_{_{-q}}(\bu){f(\bu)}\,\bd\bu, \quad\text{since}
\int_{\SQ}\chi_{_{q}}(\bu)\,\bd\bu=\delta_{{_0q}},
\end{equation}
we see the series defining the function converges on the annulus. We remark the difference in the way to compute the coefficients of a Laurent series and the usual one for functions of one complex variable.
\end{proof}

\begin{corollary}
[{\bf Holomorphic Picard Theorem for $\PQ$}]\label{holPic} Any holomorphic
line bundle $p:L\to\PQ$ determines a number $q\in\Q$ such that $L$ is holomorphically equivalent to $\I_{\C P^1_\Q}(q)$.
\end{corollary} 

 \begin{corollary}[\bf Rigidity of holomorphic vector bundles over $\PQ$] The
 holomorphic equivalence classes of holomorphic vector bundles of rank $n\in\N$ over $\PQ$ are completely determined by their
 topological invariants, i.e, the Chern classes $(q_1,\cdots,q_n)$. So their moduli spaces are discrete.
  \end{corollary}
  
 \begin{example}[\bf Vector bundles over $\C P^{n}_\Q$] Examples of holomorphic vector bundles over the $n$-dimensional adelic projective space $\C P^{n}_\Q$ can be obtained via the pullback of vector bundles
 over $\C P^n$ by the solenoidization map $p:\C P^{n}_\Q\to\C P^n$. For instance
 the moduli space of stable holomorphic vector vector bundles over 
 $\C P^2$ which decompose as the sum of two line bundles of Chern classes 
 $0$ and $m\geq2$ has a smooth, connected and rational moduli space $M(0,m)$ with 
 $\dim_\C{M(0,m)}=4m-3$ \cite{B}. By pullback we have that the moduli space
 of rank 2 holomorphic vector bundles over $\C P^{2}_\Q$ is in general non-discrete
 An interesting example for $\C P^4_\Q$ is obtained by pulling-back of the
  of the Horrocks-Mumford rank-2 vector bundle over $\C P^4$ \cite{HMu}.
 \end{example}
  
\section{Concluding remarks}\label{concrmks} 
We did not touch the very interesting topics about the adelic infinite Grassmannians and their Plücker embeddings in infinite projective spaces,  Bruhat decompositions, the unitary representations of the adelic loop groups, the 
Borel–Weil–Bott type theorems and the central extensions as well as Quillen/Segal determinant line bundles of the adelic loop groups. Also, it would be very interesting to develop 
{\it sheaf theory} for the adelic varieties studied in this article.
In particular, to develop holomorphic vector bundles over adelic toric varieties
of higher dimension. One could also try to prove theorems like Riemann-Roch for these varieties, etc. 
Bundles on $\PQ$ given by cocycles $f:A(r,R)\to\gln$ such
that the elements of the matrix $f(z)$ is given by a Laurent-Puiseux
polynomial (definition (\ref{LPPoly})) should be considered as {\it algebraic}
vector bundles. Is there a GAGA theory {\it à la Serre}? It would also be very interesting to define the canonical bundle of adelic toric varieties and
develop the notion of Calabi-Yau varieties in this adelic category. What would
be the notion of {\it mirror symmetry}?
Many other ideas come to mind.

\medskip

\centerline{\scshape Juan M. Burgos}
\medskip
{\footnotesize
 \centerline{Departamento de Matem\'aticas, CINVESTAV--CONACYT.}
   \centerline{Av. Instituto Polit\'ecnico Nacional 2508, Col. San Pedro Zacatenco, 07360, Mexico City}
  \centerline{\email{burgos@math.cinvestav.mx}}}

\medskip

\centerline{\scshape Alberto Verjovsky}
\medskip
{\footnotesize

 \centerline{Instituto de Matem\'{a}ticas} 
 \centerline{Unidad Cuernavaca, Universidad Nacional Aut\'{o}noma de M\'{e}xico,}
   \centerline{Av. Universidad S/N, C.P. 62210, Cuernavaca, Morelos, M\'{e}xico.}
}
\centerline{\email{alberto@matcuer.unam.mx}}}

\end{document}